\documentclass[english,twocolumn,final]{IEEEtran}
\usepackage[LGR,T1]{fontenc}
\usepackage[latin9]{inputenc}
\usepackage{xcolor}
\usepackage{float}
\usepackage{amsmath}
\usepackage{amssymb}
\usepackage{graphicx}
\PassOptionsToPackage{normalem}{ulem}
\usepackage{ulem}
\usepackage{epstopdf}
\usepackage{epsfig}
\usepackage{epsf}
\usepackage{graphics}

\makeatletter

\DeclareRobustCommand{\greektext}{%
  \fontencoding{LGR}\selectfont\def\encodingdefault{LGR}}
\DeclareRobustCommand{\textgreek}[1]{\leavevmode{\greektext #1}}
\ProvideTextCommand{\~}{LGR}[1]{\char126#1}

\floatstyle{ruled}
\newfloat{algorithm}{tbp}{loa}
\providecommand{\algorithmname}{Algorithm}
\floatname{algorithm}{\protect\algorithmname}
\providecolor{lyxadded}{rgb}{0,0,1}
\providecolor{lyxdeleted}{rgb}{1,0,0}

\DeclareRobustCommand{\lyxsout}[1]{\ifx\\#1\else\sout{#1}\fi}

\usepackage{bbm}

\renewcommand{\t}{^{\mbox{\tiny {T}}}}

\newcommand{\eproof}{\hfill\rule{2mm}{2mm}}

\newcommand{\bstate}{\medskip\begin{state} }
\newcommand{\estate}{ \hfill  \rule{1mm}{2mm}\medskip\end{state}}

\newcommand{\bass}{\medskip\begin{ass} }
\newcommand{\eass}{ \hfill  \rule{1mm}{2mm}\medskip\end{ass}}

\newcommand{\brem}{\medskip \begin{remark}  }
\newcommand{\erem}{\hfill \rule{1mm}{2mm}\medskip
\end{remark} }
\newcommand{\bthm}{\medskip\begin{theorem}  }
\newcommand{\ethm}{ \hfill  \rule{1mm}{2mm} \medskip
\end{theorem} }
\newcommand{\blem}{\medskip\begin{lemma}  }
\newcommand{\elem}{ \hfill \rule{1mm}{2mm}\medskip
\end{lemma} }
\newcommand{\bcorollary}{\medskip\begin{corollary}  }
\newcommand{\ecorollary}{  \hfill \rule{1mm}{2mm}\medskip
\end{corollary} }
\newcommand{\bdefn}{\medskip\begin{definition}}
\newcommand{\edefn}{  \hfill \rule{1mm}{2mm}\medskip
\end{definition} }
\newcommand{\bproposition}{\medskip\begin{proposition} }
\newcommand{\eproposition}{\hfill \rule{1mm}{2mm}\medskip
\end{proposition} }
\newcommand{\bexample}{\medskip\begin{example} \rm}
\newcommand{\eexample}{ \hfill \rule{1mm}{2mm}\medskip
\end{example} }
\newcommand{\proofnow}{\noindent{\bf Proof: }}
\newcommand{\prooflater}[1]{\noindent{\bf Proof of #1: }}

\newtheorem{theorem}{\bf Theorem}[section]
\newtheorem{ass}{\bf Assumption}[section]
\newtheorem{lemma}{\bf Lemma}[section]
\newtheorem{definition}{\bf Definition}[section]
\newtheorem{remark}{\bf Remark}[section]
\newtheorem{corollary}{\bf Corollary}[section]
\newtheorem{proposition}{\bf Proposition}[section]
\newtheorem{example}{\bf Example}[section]
\newtheorem{state}{\bf Assumption}[section]

\DeclareFontFamily{OMX}{yhex}{}
\DeclareFontShape{OMX}{yhex}{m}{n}{<->yhcmex10}{}
\DeclareSymbolFont{yhlargesymbols}{OMX}{yhex}{m}{n}
\DeclareMathAccent{\wideparen}{\mathord}{yhlargesymbols}{"F3}

\makeatother

\usepackage{babel}
\begin{document}
\title{A Sampling Control Framework and Applications to Robust and Adaptive
Control}

\author{Lijun~Zhu, and Zhiyong~Chen,~\IEEEmembership{Senior Member,~IEEE}
\thanks{ Lijun Zhu is with Key Laboratory of Imaging Processing and Intelligence Control (The Ministry of Education), and also with School of Artificial Intelligence and Automation, Huazhong University of Science and Technology, Wuhan 430074, China (e-mail:{ljzhu}@hust.edu.cn). Zhiyong Chen is with the School of Engineering, The University of Newcastle, Callaghan, NSW 2308, Australia (e-mail:zhiyong.chen@newcastle.edu.au).
%
}
}
\maketitle
\begin{abstract}
In this paper, we propose a novel sampling control framework based
on the emulation technique where the sampling error is regarded as an
auxiliary input to the emulated system. Utilizing the supremum norm of sampling
error, the design of periodic sampling
and event-triggered control law  renders the error dynamics bounded-input-bounded-state (BIBS),
and when coupled with system dynamics, achieves global or semi-global
stabilization. The proposed framework is then extended to tackle the
event-triggered and periodic sampling stabilization for  a system where
only partial state is available for feedback and the system is
subject to parameter uncertainties. The proposed framework is further
extended to solve two classes of event-triggered adaptive control
problems where the emulated closed-loop system does not admit an  input-to-state stability (ISS)
Lyapunov function. For the first class of systems with linear parameterized
uncertainties, even-triggered global adaptive stabilization is achieved
without the global Lipschitz condition on nonlinearities as often
required in the literature. For the second class of systems with uncertainties
whose bound is unknown, the event-triggered adaptive (dynamic) gain
 controller is designed for the first time. Finally, theoretical results
are verified by two numerical examples.
\end{abstract}

\section{Introduction}

 The majority of modern control relies on its digital implementation
in microprocessors and/or is deployed in a networked environment. 
And the sampled-data control scheme \cite{araki1993recent,chen2012optimal}
thus arises from the demand of more efficient controller execution
to reduce computation cost and save communication bandwidth in applications,
such as multi-robot systems, electrical power systems and chemical
processes. Sampled-data control schedules the update of a digital
controller output at periodic or aperiodic sampling instances. Input-delay
approach \cite{fridman2005input,Fridman2004,Fridman08a} is usually
adopted to design sampled-data controller for linear systems, while
two other methods are mainly used for nonlinear systems, namely controller
emulation, involving digital implementation of a continuous-time stabilizing
control law, and plant discretization consists of discretizing the
plant model and discrete-time control law design \cite{moheimani2001perspectives}.
However, the sampled-data control approach may over-sample in some
cases and motivates event-triggered control which suggests scheduling
update only when necessary and specified by the occurrence of the
designated triggering event and may achieve more efficient sampling
pattern. Sampled-data and event-triggered have been developed for
stabilization and tracking of individual systems, e.g., \cite{Tallapragada2013,Marchand2013,Tabuada2007,Liu2015}
and cooperative control of networked systems, e.g., \cite{Fan2013,Seyboth2013,chen2020often,yi2018dynamic}.
In this paper, these two control methods are uniformly called the
sampling control. 

The two-step digital emulation, constituting continuous-time control
law design and its digital implementation, is a common technique for
analysis and design of sampling control systems especially for nonlinear
systems. In \cite{Karafyllis2009,Nesic2009}, when the emulation of
periodic sampled-data control is formulated in the hybrid system setting,
the maximum allowable sampling period (MASP) that guarantees asymptotic
stability of sampled-data systems can be explicitly computed. In \cite{Chen2014a},
a new small gain theorem is proposed to reveal the quantitative tradeoff
between robustness and sampling bandwidth of the sampled-data control.
Emulation is also commonly adopted for the design of event-triggered
laws where the continuous-time controller is assumed to ensure dissipativity
\cite{moheimani2001perspectives} or the input-to-state stability
(ISS) of the closed-loop system with the sampling error as input,
see, e.g., \cite{Anta2010,DePersis2011,Liu2015,Mazo2010,Liu2015a,Tabuada2007}.
The ISS condition can be specified in a max-form, e.g. \cite{Liu2015,zhu2021}
or in an ISS-Lyapunov form, e.g. \cite{Anta2010,Liu2015a,Mazo2010,Tabuada2007}
for the closed-loop sampling system and the small gain conditions
were then proposed to ensure the stability of the event-triggered
system. In \cite{liberzon2014lyapunov}, the event-triggered technique
in \cite{Tabuada2007} was interpreted as a stabilization problem
of interconnected hybrid systems for which each subsystem admits an
ISS-Lyapunov function and a hybrid small gain condition was proposed.
As another variant of small gain theorem, the cyclic small gain theorem
has been proved effective for event-triggered control of large-scale
systems \cite{DePersis2011,Liu2015a,Liu2015}. 

Achieving Zeno free sampling and robustness to the disturbance and
uncertainties simultaneously is not trivial. In \cite{dolk2017output},
a dynamic output-based event-triggered law is proposed to achieve
a finite $L_{p}$-gain and a strictly positive lower bound on the
inter-event times when the system is subjected to non-vanishing external
disturbances. In \cite{Liu2015}, event-triggered control scheme is
designed for systems subject to disturbances. In order to handle the
uncertainties, \cite{Liu2015a} proposed a dynamic event-triggered
controller to achieve the stabilization when the system has the dynamic
uncertainties whose state is not available for the feedback. The problem
is solved in \cite{zhu2021} by a static event-triggered controller
using a particular small-gain theorem method.

Recently, a few works on the event-triggered adaptive control have
been proposed when  systems possess parameter uncertainties. Since
ISS condition is not guaranteed by the continuous-time stabilization
controller and how to design an event-triggered adaptive control scheme
is still challenging problem. In \cite{xing2016}, an event-triggered
adaptive control scheme was proposed for a class of nonlinear systems
based on the $\sigma$-modification scheme. $\sigma$-modification
scheme results in an ISS Lyapunov function, but only practical stabilization
can be achieved, i.e., the trajectories converge to a bounded set
in the neighborhood of the origin. The event-triggered law in \cite{xing2016}
was improved in \cite{huang2019} in the sense that the adaptation
dynamics was also sampled on the assumption of global Lipschitz conditions.
In \cite{wang2019}, a novel event-triggered adaptive control scheme
is proposed for a class of nonlinear systems with unknown control
direction and unknown sensor faults, while only practical stabilization
is achieved. 

In this paper, we propose a novel sampling control framework that
can lead to event-triggered and periodic sampling robust and adaptive
stabilization of nonlinear systems. The main contribution of this
paper is three-fold. First, the paper proposes a sampling control
framework where the error caused by the sampling of the actuation
is regarded as the auxiliary input to the emulated system. The emulated
closed-loop system is first assumed to admit an ISS Lyapunov function
when regarding the auxiliary input as the input, while the error dynamics
admits an continuously differentiable function whose derivative is bounded
by functions of state and sampling errors. The design of periodic
sampling and event-triggered control law renders the error dynamics
bounded-input-bounded-state (BIBS), when coupled with system dynamics,
achieves global or semi-global stabilization. Zeno free behavior is
theoretically guaranteed and the system performance is robust to the
external disturbances. Second, the proposed framework is then extended
to tackle the event-triggered and periodic sampling control for the
systems where only partial state is available for the feedback and
the system is subject to parameter uncertainties. The proposed controller
design combines the method proposed in Section \ref{sec:SCD} and
changing supply function method. In \cite{Liu2015,Liu2015a}, an auxiliary
dynamic system is introduced to estimate the decay rate of immeasurable
states and used as the dynamic threshold for the event-triggered law,
while the proposed controller is static and easier to design and implement
in practice. The periodic sampling control is also proposed here,
which is not considered in \cite{Liu2015,Liu2015a,zhu2021}. Third,
the proposed framework is extended to solve two classes of event-triggered
adaptive control problems where the emulated closed-loop system does
not admit an ISS Lyapunov function. The first class of systems contains
linear parameterized uncertainties, and global asymptotical stabilization
is achieved without the global Lipschitz condition as often required
in \cite{huang2019}. The second class of systems has the uncertainties
whose bound is not known and requires adaptive/dynamic gain technique.
To author's best knowledge, the event-triggered dynamic gain control
is first time solved in this paper.

The rest of the paper is organized as follows. Section \ref{sec:PF}
formulates the sampling control problem for nonlinear systems which
will be solved by the sampling control framework proposed in Section
\ref{sec:SCD} when an ISS condition is assumed. The proposed framework
is extended to tackle the event-triggered and periodic sampling control
for the systems where only partial state is available for the feedback
in Section \ref{sec:robust}. It is further extended to solve event-triggered
adaptive stabilization of two types of uncertain systems in Section
\ref{sec:adaptive}. The numerical simulation is conducted in Section
\ref{sec:ne} and the paper is concluded in Section \ref{sec:con}.

\textbf{Notations.} $\mathbb{R}_{\geq0}$ and $\mathbb{R}_{>0}$ denote the set of non-negative and positive real numbers, respectively. $\mathbb{R}^{n}$ denotes the real space of dimension
$n$, $\mathbb{R}^{n\times m}$ denotes the set  of real matrices with dimension
$n\times m$,  and
$\mathbb{N}$ denotes the set of  non-negative integers. Let $\|x_{[a,b]}\|:=\sup_{t\in[a,b]}\|x(t)\|$
be the supremum norm of a given signal $x:\mathbb{R}_{\geq0}\mapsto\mathbb{R}^{n}$,
over the interval $[a,b]$, and $x_{[a,b]}$ be the signal slice of
$x$ over the interval $[a,b]$. Denote a function by $\alpha\in\mathcal{K}_{\infty}$
if it is a $\mathcal{K}_{\infty}$ function, by $\alpha\in\mathcal{KL}$
if it is a $\mathcal{KL}$ function, by $m\in\mathcal{SP}$  if it
is a smooth positive function and by $m\in\mathcal{SN}$  if it
is a smooth non-negative function. 
A continuous function $\alpha:[0,\infty)\times(0,\infty)\mapsto [0,\infty)$
is said to be a parameterized $\mathcal{K}_{\infty}$ function if,
for each fixed $s>0$, the function $\alpha(\cdot,s)$ is a $\mathcal{K}_{\infty}$
function, and for each fixed $r>0$, the function $\alpha(r,\cdot)$
is non decreasing. 
For $\mathcal{K}_{\infty}$ functions $\alpha'(s)$
and $\alpha(s)$, $\alpha'(s)=\mathcal{O}[\alpha(s)]$ as $s\rightarrow0^{+}$
means that $\lim\sup_{s\rightarrow0^{+}}[\alpha'(s)/\alpha(s)]<\infty$.
Let $V:\mathbb{R}^{n}\mapsto\mathbb{R}_{\geq0}$ be a continuously
differentiable function. It is called an ISS-Lyapunov function for
the system $\dot{x}=f(x,u,d)$ if the derivative of $V(x)$ along
the $x$-dynamics satisfies, for all $x\in\mathbb{R}^{n}$, $\varpi\in\mathbb{R}^{m}$
and $d\in\mathbb{D}$,
\begin{gather}
\underline{\alpha}(\|x\|)\leq V(x)\leq\bar{\alpha}(\|x\|),\nonumber \\
\dot{V}(x)\leq-\alpha(\|x\|)+\sigma(\|\varpi\|)\label{eq:ISS-Lyap}
\end{gather}
where $\underline{\alpha},\bar{\alpha},\alpha,\sigma\in\mathcal{K}_{\infty}$.
The inequalities in  (\ref{eq:ISS-Lyap}) are simplified as $V(x)\sim\{\underline{\alpha},\bar{\alpha},\alpha,\sigma\mid\dot{x}=f(x,u,d)\}$.
A bounded piecewise continuous function $f:[0,\infty)\mapsto\mathbb{R}^{n}$
is said to be persistent exciting (PE) if there exist positive constant
$\epsilon$, $t_{0}$, $T_{0}$ such that, $\frac{1}{T_{0}}\int_{t}^{t+T_{0}}f(s)f\t(s)ds\geq\epsilon^{2}I,\;\forall t\geq t_{0}.$ 

\section{Problem Formulation of Sampling Control \label{sec:PF}}

Consider the nonlinear system 
\begin{eqnarray}
\dot{x}(t) & =& f(x(t),u(t),d(t)),\nonumber \\
y(t) & = & h(x(t)),\label{eq:dyn}
\end{eqnarray}
where $x\in\mathbb{R}^{n}$ is the system state, $y\in\mathbb{R}^{q}$
the output, $u\in\mathbb{R}^{m}$ the control input, and $d\in\mathbb{R}^{n_{d}}$
system uncertainty. Suppose the function $f:\mathbb{R}^{n}\times\mathbb{R}^{m}\times\mathbb{R}^{n_{d}}\mapsto\mathbb{R}^{n}$
is a continuously differentiable function satisfying $f(0,0,d(t))=0$
and $h(0)=0$, and the uncertainty $d$ belongs to a compact set $d\in\mathbb{D}$.
Suppose the equilibrium point of the system can be stabilized by a
continuous-time feedback controller
\begin{equation}
u(t)=\kappa(y(t)),\label{eq:ct_con}
\end{equation}
with a continuously differentiable function $\kappa$. It becomes a continuous-time
state feedback system when $y=x$. In this paper, we use the emulation
technique to study the sampling version of (\ref{eq:ct_con})
as follows
\begin{align}
u(t) & =\kappa(y(t_{k})),\;t\in[t_{k},t_{k+1}),\;k\in\mathbb{N}\label{eq:sd_con_gen}
\end{align}
where $\{t_{k}\}_{k\in\mathbb{N}}$ is a sequence of sampling
time instances. The sampling instances can be determined by an event-triggered law specified as follows 
\begin{equation}
t_{k+1}=\inf_{t>t_{k}} \{\Xi (t,t_{k},x_{[t_{k},t]}) \geq 0 \}
\label{eq:sampling-law}
\end{equation}
for some functional $\Xi$ representing event occurrence to be designed. A periodic sampling
law with sampling period $T>0$ can be considered as a special case of  
(\ref{eq:sampling-law}) with $\Xi (t,t_{k},x_{[t_{k},t]}) =t-t_k -T$, i.e., 
$t_{k+1}=t_{k}+T$.  Note
that the error caused by the sampling mechanism is
\begin{equation}
\varpi(t):=\kappa\circ h(x(t_{k}))-\kappa\circ h(x(t)),\;t\in[t_{k},t_{k+1}).\label{eq:varpi}
\end{equation}
For the emulation, the closed-loop system composed of the dynamics
(\ref{eq:dyn}) and a sampling version of feedback controller
(\ref{eq:sd_con_gen}) is put into the impulsive form 
\begin{gather}
\dot{x}=f_{c}(x,\varpi,d):=f(x,\kappa\circ h(x)+\varpi,d),\nonumber \\
\dot{\varpi}=\psi(x,\varpi,d):=-\frac{\partial \kappa\circ h(x) }{\partial x}f_{c}(x,\varpi,d),\forall t\in[t_{k},t_{k+1}),\nonumber \\
\varpi(t_{k}^{+})=0.\label{eq:ET_dyn}
\end{gather}

The objective of the sampling control is to design the sampling law
(\ref{eq:sampling-law}), either periodic sampling or event-triggered
law, such that the closed-loop system composed of (\ref{eq:dyn})
and (\ref{eq:sd_con_gen}) achieves two objectives,
\begin{enumerate}
\item stabilization: the original system is asymptotically stable at the
origin, in particular,  $\lim_{t\rightarrow\infty}x(t)=0$, semi-globally or
globally; 
\item Zeno-free behavior:  a finite number of events are triggered in a finite amount of time.
\end{enumerate}

{ 
	In general, the continuous-time feedback controller in (\ref{eq:ct_con})
	is assumed to ensure that the $x$-subsystem in (\ref{eq:ET_dyn})
	has some properties such as input-to-state stability (ISS) when regarding
	sampling error $\varpi$ as the input. 
The motivation of using
the sampling error of the actuation $\varpi $  rather than state or output measurement    is explained in \cite{khan2019}.  It is exploited in our previous papers \cite{khan2019,zhu2021}
 to design event-triggered control law 
where $x$-dynamics is assumed to have the ISS condition  specified in a max-form. As is known, the classic
  analysis tools developed for the  adaptive control  rely on the 
Lyapunov function method, typically  along the gradient of a Lyapunov function. In order to solve  event-triggered adaptive control in Section \ref{sec:adaptive}, this paper first considers  the ISS  condition specified in terms of ISS-Lyapunov function and assumes the existence of a positive function  on $\varpi$-subsystem  which plays a similar role of a Lyapunov function.
 As a result, the sampling control framework proposed in next section can be naturally extended to event-triggered adaptive control in Section \ref{sec:adaptive} where the continuous-time adaptive control  does not admit an ISS condition.
}

\bass \label{ass:SP-Bound} The $x$-subsystem in (\ref{eq:ET_dyn})
has an ISS-Lyapunov function $V(x)\sim\{\underline{\alpha},\bar{\alpha},\alpha,\sigma\mid\dot{x}=f_{c}(x,\varpi,d)\}$.
There exists a continuously differentiable $U(\varpi):\mathbb{R}^{m}\mapsto\mathbb{R}_{\geq0}$
such that the derivative of $U(\varpi)$ along the $\varpi$-subsystem,
for all $x\in\mathbb{R}^{n}$, $\varpi\in\mathbb{R}^{m}$ and $d\in\mathbb{D}$,
satisfies
\begin{gather}
\underline{\alpha}_{\varpi}(\|\varpi\|)\leq U(\varpi)\leq\bar{\alpha}_{\varpi}(\|\varpi\|),\nonumber \\
\dot{U}(\varpi)=\left(\frac{\partial U(\varpi)}{\partial\varpi}\right)\t\psi(x,\varpi,d)\leq\alpha_{\varpi}(\|\varpi\|)+\sigma_{\varpi}(\|x\|),\label{eq:bcon}
\end{gather}
for some functions $\underline{\alpha}_{\varpi}$, $\bar{\alpha}_{\varpi}$,
$\alpha_{\varpi}$, $\sigma_{\varpi}$ $\in\mathcal{K}_{\infty}$.
\eass 

\brem \label{rem:bcon} Note that the existence of the function $U(\varpi)$
condition is mild. Due to $\kappa(0)=0$ and thus $f_{c}(0,0,d)=0$,
the function $\psi$ in (\ref{eq:ET_dyn}) is continuously differentiable
and satisfies $\psi(0,0,d)=0$. Applying Lemma 11.1 in \cite{Chen2015Book},
one has $\|\psi(x,\varpi,d)\|\leq m_{1}(x,d)\|x\|+m_{2}(\varpi,d)\|\varpi\|$
for some functions $m_{1}$, $m_{2}\in\mathcal{SN}$. As a result,
$\dot{U}(\varpi)\leq m_{1}(x,d)\|\partial U(\varpi)/\partial\varpi\|\|x\|+m_{2}(\varpi,d)\|\partial U(\varpi)/\partial\varpi\|\|\varpi\|.$
Then, it is always possible to find functions $\alpha_{\varpi}$,
$\sigma_{\varpi}$ $\in\mathcal{K}_{\infty}$ such that the second inequality of
(\ref{eq:bcon}) holds. \erem

\brem
In \cite{zhu2021},  a pair of ISS and  input-to-output stability (IOS)  conditions are given in max-norm form for the closed-loop system (\ref{eq:ET_dyn}) where   $\varpi$ is regarded as the input and
$\psi$  as the output.  ISS and  IOS properties, for any $x(t_{0})$, are stated as follows,
\begin{align}
\|x(t)\| & \leq  \max\{\tilde{\beta}(\|x(t_{0})\|,t-t_{0}),\tilde{\gamma}(\|\varpi_{[t_{0},t]}\|)\},\label{eq:xISS}\\
\|\psi(t)\| & \leq  \max\{\beta(\|x(t_{0})\|,t-t_{0}),\gamma(\|\varpi_{[t_{0},t]}\|)\},\label{eq:xiIOS}
\end{align}
for $t\geq t_{0}$, where  $\gamma,\tilde{\gamma}\in\mathcal{K}_{\infty}$ and  $\beta,\tilde{\beta}$  are $\mathcal{KL}$ functions. Although ISS-Lyapunov function $V(x)$ in Assumption \ref{ass:SP-Bound} implies  (\ref{eq:xISS}), (\ref{eq:xiIOS}) and (\ref{eq:bcon}) do not necessarily imply  each other. Therefore,   conditions in Assumption \ref{ass:SP-Bound} are different from  (\ref{eq:xISS}) and  (\ref{eq:xiIOS}).   
\erem

\section{Sampling Control Design \label{sec:SCD} }
In this section, periodic sampling and event-triggered control are
proposed for the emulation system (\ref{eq:ET_dyn}). We first present
Theorem \ref{lem:uniform_lem} for the design of event-triggered law,
which provides a common guideline for periodic sampling control design. 

\subsection{Event-triggered Control}
In order to solve the sampling stabilization problem, we apply the changing
supply function technique (\cite{Sontag1995} and Lemma 2.5 of \cite{Chen2015Book})
to the ISS-Lyapunov function $V(x)$ in Assumption~\ref{ass:SP-Bound}.
Let $\mathcal{K}_\infty$ function $\alpha_{q}(s)$ be selected such that $\alpha_{q}(s)>2\sigma_{\varpi}(s)$. 
If $\alpha_{q}(s)=\mathcal{O}(\alpha(s))$ as $s\rightarrow0^{+}$,
the changing supply function technique shows that there exists another
ISS Lyapunov function $V_{q}(x)\sim\{\underline{\alpha}_{q},\bar{\alpha}_{q},\alpha_{q},\sigma_{q}\mid\dot{x}=f_{c}(x,\varpi,d)\}$
where the functions $\underline{\alpha}_{q}$, $\bar{\alpha}_{q},\sigma_{q}\in\mathcal{K}_{\infty}$
can be calculated accordingly. Then, the event-triggered control law
is presented as follows.

\bthm \label{lem:uniform_lem} Suppose the system composed of (\ref{eq:dyn})
and (\ref{eq:sd_con_gen}) satisfies Assumption \ref{ass:SP-Bound}.
Let $\mathcal{K}_\infty$ functions $\hat{\alpha}_{\varpi}$ and $\gamma$ be  $\hat{\alpha}_{\varpi}(s)>\sigma_{q}(s)/2$ and $\gamma(s)>\hat{\alpha}_{\varpi}(s)+\alpha_{\varpi}(s)$.
Suppose $\varpi(t_{0})=0$, $\sigma_{\varpi}(s)=\mathcal{O}(\alpha(s))$
and $\gamma(s)=\mathcal{O}(\underline{\alpha}_{\varpi}(s))$ as $s\rightarrow0^{+}$.
If the event-triggered law in (\ref{eq:sampling-law}) is designed
as 
\begin{gather}
 t_{k+1}=\inf_{t\geq t_{k}}\{2(t-t_{k})\gamma(\|\varpi_{[t_{k},t]}\|)\nonumber \\
\geq\max_{\tau\in[t_{k},t]}\{U(\varpi(\tau))\},\;\text{and }\|\varpi_{[t_{k},t]}\|\neq0\},\label{eq:sd_con}
\end{gather}
the equilibrium
point $x=0$ of the system is globally asymptotically stable and  Zeno-behavior is avoided. Moreover, $V_{q}(x(t))\leq V_{q}(x(t_{k}))$,
$\forall t\in[t_{k},t_{k+1}]$ and $k\in\mathbb{N}$. 
\ethm


\proofnow The closed-loop system consisting of (\ref{eq:dyn}) and
(\ref{eq:sd_con_gen}) is put into the form (\ref{eq:ET_dyn}).
The proof will be divided into three parts. First, the boundedness
of all signals is proved. The fact that $\sigma_{\varpi}(s)=\mathcal{O}(\alpha(s))$
as $s\rightarrow0^{+}$ implies that it is always possible to find a new
supply function $V_{q}(x)\sim\{\underline{\alpha}_{q},\bar{\alpha}_{q},\alpha_{q},\sigma_{q}\mid\dot{x}=f_{c}(x,\varpi,d)\}$
and $\alpha_{q}(s)>2\sigma_{\varpi}(s)$. Integrating both side of
(\ref{eq:bcon}) in Assumption \ref{ass:SP-Bound} gives
\begin{gather*}
U(\varpi(t))=\int_{t_{k}}^{t}\dot{U}(\varpi(\tau))d\tau\leq\int_{t_{k}}^{t}\alpha_{\varpi}(\|\varpi(\tau)\|)d\tau\\
+\int_{t_{k}}^{t}\sigma_{\varpi}(\|x(\tau)\|)d\tau,\;t\in[t_{k},t_{k+1}]
\end{gather*}
which leads to
\begin{gather}
\max_{\tau\in[t_{k},t]}\{U(\varpi(\tau))\}\leq\int_{t_{k}}^{t}\alpha_{\varpi}(\|\varpi(\tau)\|)d\tau\nonumber \\
+\int_{t_{k}}^{t}\sigma_{\varpi}(\|x(\tau)\|)d\tau,\;t\in[t_{k},t_{k+1}]\label{eq:UW}
\end{gather}
due to $\alpha_{\varpi}(\|\varpi(\tau)\|)\geq0$ and $\sigma_{\varpi}(\|x(\tau)\|)\geq0,\forall\tau\in[t_{k},t]$.
It further implies that 
\begin{gather}
\max_{\tau\in[t_{k},t]}\{U(\varpi(\tau))\}\leq-\int_{t_{k}}^{t}\hat{\alpha}_{\varpi}(\|\varpi(\tau)\|)d\tau\nonumber \\
+\int_{t_{k}}^{t}\sigma_{\varpi}(\|x(\tau)\|)d\tau+\int_{t_{k}}^{t}\left[\alpha_{\varpi}(\|\varpi(\tau)\|)+\hat{\alpha}_{\varpi}(\|\varpi(\tau)\|)\right]d\tau\nonumber \\
\leq2\max\left\{ -\int_{t_{k}}^{t}\hat{\alpha}_{\varpi}(\|\varpi(\tau)\|)d\tau+\int_{t_{k}}^{t}\sigma_{\varpi}(\|x(\tau)\|)d\tau,\right.\nonumber \\
\left.\int_{t_{k}}^{t}\left[\alpha_{\varpi}(\|\varpi(\tau)\|)+\hat{\alpha}_{\varpi}(\|\varpi(\tau)\|)\right]d\tau\right\} ,\;\forall t\in[t_{k},t_{k+1}]\label{eq:U_ine}
\end{gather}
where $\hat{\alpha}_{\varpi}(s)>\sigma_{q}(s)/2$. The sampling law
(\ref{eq:sd_con}) implies $\int_{t_{k}}^{t}\gamma(\varpi(\tau))d\tau\leq(t-t_{k})\gamma(\|\varpi_{[t_{k},t]}\|)\leq\frac{1}{2}\max_{\tau\in[t_{k},t]}\{U(\varpi(\tau))\}$,
which together with $\gamma(s)>\hat{\alpha}_{\varpi}(s)+\alpha_{\varpi}(s)$,
leads to 
\begin{gather}
U(\varpi(t))\leq\max_{\tau\in[t_{k},t]}\{U(\varpi(\tau))\}\leq-2\int_{t_{k}}^{t}\hat{\alpha}_{\varpi}(\|\varpi(\tau)\|)d\tau\nonumber \\
+2\int_{t_{k}}^{t}\sigma_{\varpi}(\|x(\tau)\|)d\tau,\;\forall t\in[t_{k},t_{k+1}].\label{eq:Uw_ine_bound}
\end{gather}
Let $\bar{V}(x,\varpi)=V_{q}(x)+U(\varpi)$ be the Lyapunov function
candidate for the closed-loop system (\ref{eq:ET_dyn}). Since $x(t)$
and hence $V_{q}(x(t))$ are continuous in $t$ and $U(\varpi(t))$ is piecewise
continuous in $t$ and has jump at $t_{k}$, one has
\begin{gather*}
\bar{V}(x(t_{k}^{+}),\varpi(t_{k}^{+}))=V_{q}(x(t_{k}))\leq\bar{V}(x(t_{k}),\varpi(t_{k}))\\
=U(\varpi(t_{k}))+V_{q}(x(t_{k})),\;\forall k\in\mathbb{N}
\end{gather*}
by noting $\varpi(t_{k}^{+})=0$. For $t\in[t_{k}^+,t_{k+1}]$, 
\begin{gather}
\bar{V}(x(t),\varpi(t))-\bar{V}(x(t_{k}^+),\varpi(t_{k}^+))\nonumber \\
=U(\varpi(t))+V_{q}(x(t))-V_{q}(x(t_{k}^+))\nonumber\\=U(\varpi(t))+\int_{t_{k}^+}^{t}\dot{\bar{V}}(x(\tau),\varpi(\tau))d\tau\nonumber \\
\leq-2\int_{t_{k}^+}^{t}\hat{\alpha}_{\varpi}(\|\varpi(\tau)\|)d\tau+2\int_{t_{k}^+}^{t}\sigma_{\varpi}(\|x(\tau)\|)d\tau\nonumber \\
-\int_{t_{k}^+}^{t}\alpha_{q}(\|x(\tau)\|)d\tau+\int_{t_{k}^+}^{t}\sigma_{q}(\|\varpi(\tau)\|)d\tau\nonumber \\
\leq-\int_{t_{k}^+}^{t}\alpha_{\gamma}(\|x(\tau)\|)d\tau-\int_{t_{k}^+}^{t}\gamma_{\alpha}(\|\varpi(\tau)\|)d\tau,\label{eq:V_bar_ine}
\end{gather}
where the functions $\alpha_{\gamma}$ and $\gamma_{\alpha}$ are defined
as 
\begin{gather*}
\gamma_{\alpha}(s)=2\hat{\alpha}_{\varpi}(s)-\sigma_{q}(s)>0,\\
\alpha_{\gamma}(s)=\alpha_{q}(s)-2\sigma_{\varpi}(s)>0.
\end{gather*}
As a result, $\bar{V}(x(t),\varpi(t))$ is monotonically decreasing
except when $\text{col}(x,\varpi)=0$ and therefore the signals $x$ and
$\varpi$ are bounded for $t\geq t_{0}$ whose bound depends on the
initial value $x(t_{0})$. Denote the bound $\varpi$ by $R_{\varpi}$.
The derivative of $\varpi$ denoted as $q:=\dot{\varpi}$ is also
bounded whose bound is denoted by $\|q(t)\|\leq R(x(t_{0})),\forall t>t_{0}$.
The notation $R$ will be used instead of $R(x(t_{0}))$ for the notation
simplicity. Note that $V_{q}(x(t_{k}))=\bar{V}(x(t_{k}^{+}),\varpi(t_{k}^{+}))\geq\bar{V}(x(t),\varpi(t))\geq V_{q}(x(t))$
for any $t\in[t_{k},t_{k+1}]$ and $k\in\mathbb{N}$. 

Then, it will be shown that the event-triggered law is free of Zeno
behavior by showing that $t_{k+1}-t_{k}>c$ for some constant $c>0$ possibly depending on the initial condition. Note that $\gamma(s)=\mathcal{O}(\underline{\alpha}_{\varpi}(s))$
as $s\rightarrow0^{+}$ implies that $\lim_{s\rightarrow0^{+}}\gamma(s)/\underline{\alpha}_{\varpi}(s)<\infty.$
Denote $C=\lim_{s\rightarrow0^{+}}\gamma(s)/\underline{\alpha}_{\varpi}(s)$.
Due to $\varpi(t)=-\int_{t_{k}}^{t}q(\tau)d\tau,\;\forall t\in[t_{k},t_{k+1})$,
one has $\|\varpi_{[t_{k},t]}\|\leq(t-t_{k})\|q_{[t_{k},t]}\|\leq(t-t_{k})R\leq\epsilon,\;\forall t\in[t_{k},t_{k}+\frac{\epsilon}{R}]$.
For a given $\delta$, there exists an $\epsilon$ such that
\begin{equation}
\|\varpi_{[t_{k},t]}\|<\epsilon\implies|\gamma(\|\varpi_{[t_{k},t]}\|)/\underline{\alpha}_{\varpi}(\|\varpi_{[t_{k},t]}\|)-C|<\delta.\label{eq:w_inq}
\end{equation}
Denote 
\begin{equation}
t_{d}=\min\left \{\frac{1}{2(C+\delta)},\epsilon/R\right\}\label{eq:td}
\end{equation}
As a result
\begin{gather*}
2(t-t_{k})\gamma(\|\varpi_{[t_{k},t]}\|)\leq2(t-t_{k})(C+\delta)\underline{\alpha}_{\varpi}(\|\varpi_{[t_{k},t]}\|)\\
\leq\underline{\alpha}_{\varpi}(\|\varpi_{[t_{k},t]}\|)\leq\max_{\tau\in[t_{k},t]}\{U(\varpi(\tau))\},\;\forall t\in[t_{k},t_{k}+t_{d}],
\end{gather*}
where we used (\ref{eq:w_inq}). Therefore, the sampling interval
can be selected to be lower bounded irrespective of $k$, i.e., $t_{k+1}-t_{k}\geq T_{c}=t_{d}$
for $k\in\mathbb{N}$ and the bound $t_{d}$ depends on $R(x(t_0))$ and thus the
initial condition $x(t_{0})$. 

Finally, the convergence of the signal $x$ to zero, i.e., $\lim_{t\rightarrow\infty}x(t)=0$,
is proved. If it is not true, there is a positive constant $c_{1}$
such that for every $t>0$, we can find $\bar{t}\geq t$ with $\|x(\bar{t})\|\geq c_{1}$. It is noted from (\ref{eq:V_bar_ine})
 that $\bar{V}(x(t_{k+1}),\varpi(t_{k+1}))-\bar{V}(x(t_{k}),\varpi(t_{k}))\leq-\int_{t_{k}}^{t_{k+1}}\alpha_{\gamma}(\|x(\tau)\|)d\tau$
and $\lim_{t\rightarrow\infty}\int_{t_{0}}^{t}\alpha_{\gamma}(\|x(\tau)\|)d\tau$
exists and is finite. Since $x(t)$ is a uniformly continuous function,
applying Barbalat\textquoteright s lemma we can show that $\lim_{t\rightarrow\infty}x(t)=0$.
It will also be shown that the convergence of the signal $\varpi$
to zero, i.e., $\lim_{t\rightarrow\infty}\varpi(t)=0$. If it is not
true, there is a positive constant $c_{2}$ such that for every $t>0$,
we can find $\bar{t}\geq t$ with $\|\varpi(\bar{t})\|\geq c_{2}$.
Also, there exists a $k=\max_{k}\{k|t_{k}<\bar{t}\text{ and }\varpi(t_{k}^+)=0\}$.
Since signals $x$ and $\varpi$ are bounded and $\varpi$ is continuous
in $t\in[t_{k},t_{k+1}]$, there exists a $\tilde{t}\in(t_{k},\bar{t})$
such that $\|\varpi(\tilde{t})\|=\frac{1}{2}c_{2}$, $\|\varpi(\tilde{t}_{p})\|\geq\frac{1}{2}c_{2},$
$\forall\tilde{t}_{p}\in[\tilde{t},\bar{t}]$ and $\bar{t}-\tilde{t}\geq c_{3}$
for some $c_{3}>0$. Following the proof of Barbalat\textquoteright s
lemma, it causes a contradiction and therefore $\lim_{t\rightarrow\infty}\varpi(t)=0$.
\eproof

\brem \label{rem:bound} The bound of $\varpi$ is explicitly derived
below. The inequality (\ref{eq:UW}) and the event-triggered law (\ref{eq:sd_con})
leads to 
\begin{gather*}
2(t-t_{k})\alpha_{\varpi}(\|\varpi_{[t_{k},t]}\|)\leq2(t-t_{k})\gamma(\|\varpi_{[t_{k},t]}\|)\\
\leq\max_{\tau\in[t_{k},t]}\{U(\varpi(\tau))\} \\
\leq\int_{t_{k}}^{t}\alpha_{\varpi}(\|\varpi(\tau)\|)d\tau+\int_{t_{k}}^{t}\sigma_{\varpi}(\|x(\tau)\|)d\tau\\
\leq(t-t_{k})\alpha_{\varpi}(\|\varpi\|_{[t_{k},t]})+(t-t_{k})\sigma_{\varpi}(\|x\|_{[t_{k},t]}),\;t\in[t_{k},t_{k+1}]
\end{gather*}
which further implies that
\begin{equation}
\alpha_{\varpi}(\|\varpi\|_{[t_{k},t]})\leq\sigma_{\varpi}(\|x\|_{[t_{k},t]}),\;\forall t\in[t_{k},t_{k+1}]. \label{eq:varpi_bound}
\end{equation}
Note that $\underline{\alpha}_{q}(\|x(t)\|)\leq V_{q}(x(t))\leq\bar{V}(x(t),\varpi(t))\leq\bar{V}(x(t_{0}),\varpi(t_{0}))= V_{q}(x(t_{0})),$
therefore one has $\|x(t)\|\leq\underline{\alpha}_{q}{}^{-1}\left(V_{q}(x(t_{0}))\right)$
and $\|\varpi(t)\|\leq\alpha_{\varpi}^{-1}\circ\sigma_{\varpi}\circ\underline{\alpha}_{q}{}^{-1}\left(V_{q}(x(t_{0})\right)$,
$\forall t\geq t_{0}.$  \erem

	\brem
The inequality (\ref{eq:varpi_bound})
	 implies that for any bounded input signal $x$, i.e., $\|x(t)\|\leq x_{b},\forall t\geq t_{0}$,
the trajectory of $\varpi$ stays in the
	ball $B:=\{\varpi\mid\|\varpi||\leq {\alpha}_{\varpi}^{-1}\circ\sigma_{\varpi}(x_{b})\}.$
	Therefore, the event-triggered law in (\ref{eq:sd_con}) render the
	$\varpi$-dynamics bounded-input-bounded-state when regarding $x$ as the
	input. \erem


\brem \label{rem:fun_sel} The condition $\gamma(s)=\mathcal{O}(\underline{\alpha}_{\varpi}(s))$
as $s\rightarrow0^{+}$ in Theorem \ref{lem:uniform_lem} can be implied
by two conditions $\sigma(s)=\mathcal{O}(\underline{\alpha}_{\varpi}(s))$
and $\alpha_{\varpi}(s)=\mathcal{O}(\underline{\alpha}_{\varpi}(s))$
as $s\rightarrow0^{+}$. Note that $\sigma(s)=\mathcal{O}(\underline{\alpha}_{\varpi}(s))$
as $s\rightarrow0^{+}$ implies $\sigma_{q}(s)=\mathcal{O}(\underline{\alpha}_{\varpi}(s))$
as $s\rightarrow0^{+}$. It together with  $\alpha_{\varpi}(s)=\mathcal{O}(\underline{\alpha}_{\varpi}(s))$
as $s\rightarrow0^{+}$ implies that one can find a $\gamma(s)$ such
that $\lim_{s\rightarrow0^{+}}\gamma(s)/\underline{\alpha}_{\varpi}(s)<\infty$
or equivalently $\gamma(s)=\mathcal{O}(\underline{\alpha}_{\varpi}(s))$
as $s\rightarrow0^{+}$. 

To see we can make $\alpha_{\varpi}(s)=\mathcal{O}(\underline{\alpha}_{\varpi}(s))$ as $s\rightarrow 0^+$, 
let  $U(\varpi)=1/2\|\varpi\|^{2}$ and its
derivative is calculated as 
\begin{gather}
\dot{U}(\varpi)\leq m_{1}(x,d)\|x\|\|\varpi\|+m_{2}(\varpi,d)\|\varpi\|^{2}\nonumber \\
\leq\left[\frac{1}{4}+m_{2}(\varpi,d)\right]\|\varpi\|^{2}+m_{1}^{2}(x,d)\|x\|^{2}\label{eq:Uw_ine}
\end{gather}
where we used the bound of $\|\psi(x,\varpi,d)\|$ in Remark \ref{rem:bcon}.
Then, inequality (\ref{eq:bcon}) in Assumption \ref{ass:SP-Bound}
is satisfied with $\alpha_{\varpi}(\|\varpi\|)\geq\left[\frac{1}{4}+m_{2}(\varpi,d)\right]\|\varpi\|^{2}$.
Note that in this case $\underline{\alpha}_{\varpi}(s)=1/2\|\varpi\|^{2}$
and then $\lim_{s\rightarrow0^{+}}\alpha_{\varpi}(s)/\underline{\alpha}_{\varpi}(s)>0$
is satisfied due to $m_{2}(\varpi,d)>0$. 

The  conditions $\sigma(s)=\mathcal{O}(\underline{\alpha}_{\varpi}(s))$
and $\sigma_{\varpi}(s)=\mathcal{O}(\alpha(s))$ as $s\rightarrow0^{+}$
can be made satisfied during the continuous-time feedback controller
design phase. From (\ref{eq:Uw_ine}), we can choose the function
$\sigma_{\varpi}$ in (\ref{eq:bcon}) to be $\sigma_{\varpi}(\|x\|)\geq\bar{m}(x,d)\|x\|^{2}$
for some function $\bar{m}(x,d)\geq m_{1}^{2}(x,d)>0$. Then, if one
designs the continuous-time feedback controller that renders $\lim_{s\rightarrow0^{+}}\alpha(s)/s^{2}<\infty$
and $\lim_{s\rightarrow0^{+}}\sigma(s)/s^{2}<\infty$, it leads to
$\sigma(s)=\mathcal{O}(\underline{\alpha}_{\varpi}(s))$ and $\sigma_{\varpi}(s)=\mathcal{O}(\alpha(s))$
as $s\rightarrow0^{+}$. Such kind of feedback controller can always
be found for a large class of nonlinear systems such as those in strict
feedback form and lower-triangular form, to be presented in Section
\ref{sec:robust} and \ref{sec:adaptive}. \erem

\brem \label{rem:Tk} It is observed from the proof of Theorem \ref{lem:uniform_lem}
that if the next sampling time is selected as $t_{k+1}=t_{k}+T_{k}$,
for any $T_{k}\leq\bar{t}_{k+1}-t_{k}$ where $\bar{t}_{k+1}$ is the next sampling time calculated from Theorem \ref{lem:uniform_lem}, i.e.,
\begin{gather}
\bar{t}_{k+1}:=\inf_{t\geq t_{k}}\{2(t-t_{k})\gamma(\|\varpi_{[t_{k},t]}\|)\nonumber \\
\geq\max_{\tau\in[t_{k},t]}\{U(\varpi(\tau))\},\;\text{and }\|\varpi_{[t_{k},t]}\|\neq0\}, \label{eq:t_k_bar}
\end{gather}
 results of Theorem \ref{lem:uniform_lem}
still hold.\erem

The following proposition shows that event-triggered control tends
to behave like a periodic sampling control as $t\rightarrow\infty$,
whose sampling interval approaches a constant. 

\bproposition \label{prop:et} Suppose the system is composed of
(\ref{eq:dyn}) and (\ref{eq:sd_con_gen}) and the event-triggered
control law is designed according to (\ref{eq:sd_con}) in Theorem
\ref{lem:uniform_lem}. Let $\underline{\mu}:=\lim_{s\rightarrow0^{+}}\underline{\alpha}_{\varpi}(s)/\gamma(s)$.
Then, $\lim_{k\rightarrow\infty}(t_{k+1}-t_{k})>\underline{\mu}/2$.
Moreover, if the function $U(\varpi)$ is specified as a $\mathcal{K}_{\infty}$
function, i.e., $U(\varpi)=U(\|\varpi\|)$ and $\mu:=\lim_{s\rightarrow0^{+}}U(s)/\gamma(s)$,
then, $\lim_{k\rightarrow\infty}(t_{k+1}-t_{k})=\mu/2$. \eproposition
\proofnow The proof is similar to that of Proposition 2.1 in our
paper \cite{zhu2021} and omitted here. \eproof

\subsection{Periodic Sampling Control \label{sub:PSC}}

In this subsection, the periodic sampling controller is to seek a
uniform sampling period, simply denoted as 
\begin{equation}
T:=t_{k+1}-t_{k}\label{eq:period}
\end{equation}
The constant $T$ is called the sampling period and $\omega:=1/T$
sampling frequency. From Remark \ref{rem:Tk}, looking for the 
sampling period $T$ amounts to, for all $k\in\mathbb{N}$,
finding a uniform $T$ such that $t_{k+1}=t_{k}+T\leq\bar{t}_{k+1}$ with $\bar{t}_{k+1}$ given in (\ref{eq:t_k_bar}), or equivalently,   the following inequality must be made satisfied
\begin{gather}
2(t-t_{k})\gamma(\|\varpi_{[t_{k},t]}\|)<\max_{\tau\in[t_{k},t]}\{U(\varpi(\tau))\},\nonumber \\
\forall t\in[t_{k},t_{k}+T).\label{eq:psc_con}
\end{gather}
The theorem of periodic sampling control follows the idea.
{
\bthm \label{thm:psc} Suppose the system composed of (\ref{eq:dyn})
and (\ref{eq:sd_con_gen}) satisfies Assumption \ref{ass:SP-Bound}. Let
 functions  $V_q$, $\gamma$, $\alpha_{\varpi}$, $\sigma_{\varpi}$  and $\underline{\alpha}_{q}$ 
be defined in the proof of Theorem \ref{lem:uniform_lem}. Let $\mathcal{X}\subset \mathbb{R}^n$ be a compact set and 
 $R_{0}=\max_{x\in\bar{\mathcal{X}}}\{ \alpha_{\varpi}^{-1}\circ\sigma_{\varpi}\circ\underline{\alpha}_{q}{}^{-1}\left(V_{q}(x)\right)\}$ where $\bar{\mathcal{X}}$ is the closure of $\mathcal{X}$.
 Find $T$ for (\ref{eq:period})
such that 
\begin{gather}
2T\gamma(s)<\underline{\alpha}_{\varpi}(s),\;\forall0<s\leq R_{0},\label{eq:PSC1}
\end{gather}
the equilibrium point $x=0$ of the system is asymptotically stable for any initial condition $x(t_{0})\in \mathcal{X} $. Moreover, $\mathcal{X} \subseteq \mathbb{S}$ where $\mathbb{S}=\{x|V_q(x)\leq \underline \alpha_q \circ \sigma_\varpi^{-1} \circ \alpha_\varpi(R_0)\}$ is positive invariant
and  $\|\varpi(t)\|\leq R_{0},\forall t\geq t_{0}$. \ethm }

The proof of Theorem \ref{thm:psc} is given in Appendix which will
show that the condition (\ref{eq:PSC1}) implies (\ref{eq:psc_con}).

\brem \label{rem:psc_init} For a given set $\mathcal{X}$ and thus $R_{0}$, Zeno freeness
proof in Theorem \ref{lem:uniform_lem} shows that $T$ satisfying
condition (\ref{eq:PSC1}) can always be found.  Note
that $R_{0}$ is determined by the initial state set $\mathcal{X}$ and needs to be known,
therefore the periodic sampling law in Theorem \ref{thm:psc} achieve
semi-global stabilization rather than the global stabilization. Similarly,
for a given $T>0$, one can also find $R_{0}>0$ such that the condition
(\ref{eq:PSC1}) is satisfied. Then, Theorem \ref{thm:psc} can be
interpreted in another way. That is, given a sampling limitation in
terms of the fastest sampling speed $\omega$, the best stabilization
performance can be achieved in terms of the estimated region of attraction,
i.e., $x(t_{0})\in\mathbb{S}$.
Note that the best estimated periodic sampling interval $T$ is sometimes
called maximum allowable sampling period. If the function $U(\varpi)$
can be specified as a $\mathcal{K}_{\infty}$ function, i.e., $U(\varpi)=U(\|\varpi\|)$,
then the condition can be refined as $2T\gamma(s)<U(s),\;\forall0<s\leq R_{0}$.
In this case, $T$ can be found as 
\[
T=\frac{1}{2}\min_{s\in(0,R_{0}]}\frac{U(s)}{\gamma(s)}\leq\frac{1}{2}\mu:=\lim_{s\rightarrow0^{+}}U(s)/\gamma(s)
\]
where $\mu/2$ is the asymptotic sampling period of the event-triggered
in Proposition \ref{prop:et} as $t$ approaches infinity. In a word,
as the system trajectory approaches equilibrium point, the periodic
sampling is no better than event-triggered control in terms of saving
the sampling times. \erem 

\subsection{Robustness Issue}

Now, let us consider the robustness of the sampling control when the
uncertainties $d$ also includes non-vanishing external disturbance.
In this case, we suppose the emulated system (\ref{eq:ET_dyn}) satisfies
the external-disturbance version of Assumption \ref{ass:SP-Bound}.

\bass \label{ass:SP_bound_dist} Suppose the $x$-subsystem in (\ref{eq:ET_dyn})
has an ISS-Lyapunov function $V(x)\sim\{\underline{\alpha},\bar{\alpha},\alpha,(\sigma,\varsigma)\mid\dot{x}=f_{c}(x,\varpi,d)\}$,
i.e., $\underline{\alpha}(\|x\|)\leq V(x)\leq\bar{\alpha}(\|x\|),\dot{V}(x)\leq-\alpha(\|x\|)+\sigma(\|\varpi\|)+\varsigma(\|d\|)$
where $\underline{\alpha},\bar{\alpha},\alpha,\sigma,\varsigma\in\mathcal{K}_{\infty}$.
There exists a continuously differentiable $U(\varpi):\mathbb{R}^{n}\rightarrow\mathbb{R}_{\geq0}$
such that the derivative of $U(\varpi)$ along the $\varpi$-subsystem,
for all $x\in\mathbb{R}^{n}$, $\varpi\in\mathbb{R}^{m}$ and $d\in\mathbb{D}$,
$\underline{\alpha}_{\varpi}(\|\varpi\|)\leq U(\varpi)\leq\bar{\alpha}_{\varpi}(\|\varpi\|)$
and $\dot{U}(\varpi)=\left(\frac{\partial U(\varpi)}{\partial\varpi}\right)\t\psi(x,\varpi,d)\leq\alpha_{\varpi}(\|\varpi\|)+\sigma_{\varpi}(\|x\|)+\varsigma_{\varpi}(\|d\|),$
for some $\underline{\alpha}_{\varpi}$, $\bar{\alpha}_{\varpi}$,
$\alpha_{\varpi}$, $\sigma_{\varpi},\varsigma_{\varpi}\in\mathcal{K}_{\infty}$.
\eass 

Then, we present the following corollary for event-triggered law whose proof is
given in the Appendix and the result for periodic sampling is similar and not presented.

\bcorollary \label{cor:robust} Suppose the system composed of (\ref{eq:dyn})
and (\ref{eq:sd_con_gen}). Suppose conditions of Theorem \ref{lem:uniform_lem}
are satisfied except that Assumption \ref{ass:SP-Bound} is replaced
by Assumption \ref{ass:SP_bound_dist}. Let $\bar{d}$ be the bound
of external disturbance $d(t)$, i.e., $\|d(t)\|\leq\bar{d}$ for
$t\geq t_{0}$. If the event-triggered law is designed as in Theorem
\ref{lem:uniform_lem}, then the trajectories of closed-loop system
is ultimately bounded, i.e., there exists a $\Delta t(x_{0})$ possibly depending
on the initial condition $x_{0}$ such that 
\[
\|x\|\leq\underline{\alpha}_{\upsilon}^{-1}\circ\bar{\alpha}_{\upsilon}\circ\alpha_{\upsilon}{}^{-1}\circ\varsigma_{\upsilon}(\bar{d}),\;\forall t\geq t_{0}+\Delta t(x_{0}),
\]
for some $\underline{\alpha}_{\upsilon}$, $\bar{\alpha}_{\upsilon}$,
$\alpha_{\upsilon}$, $\varsigma_{\upsilon}\in\mathcal{K}_{\infty}$ and Zeno-behavior is avoided.
\ecorollary

\section{Robust Sampling Control \label{sec:robust}}

In this section, we will apply the sampling control method proposed
in Section \ref{sec:SCD} to solve the sampling robust stabilization
problem of a class of nonlinear systems, called strict feedback systems
of a relative degree one, as follows, 
\begin{eqnarray}
\dot{z} & = & q(z,x,d)\nonumber \\
\dot{x} & = & f(z,x,d)+u\label{eq:rb_sys}
\end{eqnarray}
where $z\in\mathbb{R}^{p}$ and $x\in\mathbb{R}$ are state variables,
$u\in\mathbb{R}$ is the input and $d\in\mathbb{D}$ is the uncertainties
belonging to a compact set $\mathbb{D}\in\mathbb{R}^{l}$. The functions
$q$ and $f$ are sufficiently smooth with $q(0,0,d)=0$ and $f(0,0,d)=0$
for all $d\in\mathbb{D}$. Note that the state $z$ is assumed to
be not available for feedback control and thus the $z$-dynamics is
called dynamic uncertainty. A common assumption on $z$-dynamics is
given to make the problem tractable. 

\bass \label{ass:z} The $z$-subsystem in (\ref{eq:rb_sys}) has
an ISS-Lyapunov function $V(z)\sim\{\underline{\alpha},\bar{\alpha},\alpha,\sigma\mid\dot{z}=q(z,x,d)\}$
and functions $\alpha(s)$ and $\sigma(s)$ satisfy $\limsup_{s\rightarrow0^{+}}s^{2}/\alpha(s)<\infty,\;\limsup_{s\rightarrow0^{+}}\sigma(s)/s^{2}<\infty.$
\eass

The sampling robust stabilization problem is to design sampling controller
$u$ such that $\lim_{t\rightarrow\infty}\text{col}(z(t),x(t))=0$.
Since $f(z,x,d)$ is a sufficiently smooth function, one has 
\begin{equation}
|f(z,x,d)|\le m_{1}(z)\|z\|+m_{2}(x)|x|,\forall d\in\mathbb{D}\label{eq:f_bound}
\end{equation}
for some sufficiently smooth functions $m_{1},m_{2}\in\mathcal{SN}$
depending on the size of $\mathbb{D}$. For the continuous-time stabilization
of the system, a high-gain controller can be adopted to dominate the
uncertainties when the size of $\mathbb{D}$ is known. The case that
the size of $\mathbb{D}$ is unknown will be handled using dynamic
gain technique in Section \ref{sec:adaptive}. The continuous-time
controller usually takes the form of $u=\kappa(x)=-\rho(x)x$ with
the high-gain term $\rho(x)$ to be specified. We adopt the method
developed in Section \ref{sec:SCD} and propose the sampling controller
as follows 
\begin{align}
u(t) & =\kappa(x(t_{k}))=-\rho(x(t_{k}))x(t_{k}),\;t\in[t_{k},t_{k+1}),\label{eq:rb_con}
\end{align}
for $k\in\mathbb{N}$.
Define the sampling error $\varpi(t)$ as $\varpi(t)=\kappa(x(t_{k}),z(t_{k}))-\kappa(x(t),z(t)),\;t\in[t_{k},t_{k+1}).$
Then, the sampled-data closed-loop system is rewritten as
\begin{eqnarray}
\dot{z} & = & q(z,x,d)\nonumber \\
\dot{x} & = & f(z,x,d)-\rho(x)x+\varpi,\nonumber \\
\dot{\varpi} & = &  \frac{d \kappa(x)}{dx} (f(z,x,d)-\rho(x)x+\varpi).
\label{eq:rb_sys_cl}
\end{eqnarray}
The event-triggered robust stabilization can be solved by applying Theorem \ref{lem:uniform_lem}. 

\bthm \label{thm:robust} Consider the system composed of (\ref{eq:rb_sys})
and (\ref{eq:rb_con}) under Assumption \ref{ass:z}. There exists
sufficiently smooth positive functions $\rho:\mathbb{R}\mapsto\mathbb{R}_{>0}$
and $\gamma:\mathbb{R}_{\geq 0}\mapsto\mathbb{R}_{\geq 0}$ such that when
the event-triggered law 
 is designed as
\begin{gather}
t_{k+1}=\inf\{t>t_{k}|2(t-t_{k})\gamma(\|\varpi\|_{[t_k,t]})\nonumber \\
\geq\max_{\tau\in[t_{k},t]}\{U(\varpi(\tau))\},\;\text{and }\|\varpi_{[t_{k},t]}\|\neq0\},\label{eq:sd_con1-1-1}
\end{gather}
the equilibrium
point $\text{col}(z,x)=0$ of the system is globally asymptotically
stable and  Zeno-behavior is avoided.  \ethm 

\proofnow Denote $\xi=\text{col}(z,x)$. Let us first consider the
$\xi$-dynamics. Let $\Delta(z)=m_{1}^{2}(z)+1$. By changing supply
function technique, there exists another ISS Lyapunov function for
$z$-dynamics, $V_{z}(z)\sim\{\underline{\alpha}_{z},\bar{\alpha}_{z},\Delta(z)\|z\|^{2},\varkappa(x)x^{2}\mid\dot{z}=q(z,x,d)\}$
for some functions $\underline{\alpha}_{z}$, $\bar{\alpha}_{z}\in\mathcal{K}_{\infty}$,
and $\varkappa\in\mathcal{SN}$, that can be calculated accordingly.
Let
\[
\rho(x)\geq[\varkappa(x)+m_{2}(x)+3/2]
\]
and the Lyapunov function candidate be $V_{\xi}(\xi)=V_{z}(z)+x^{2}/2$.
The calculation of the derivative of $V_{\xi}(\xi)$, along the trajectory
of $\xi$-dynamics, obtains
\begin{gather*}
\dot{V}_{\xi}(\xi)\leq-\Delta(z)\|z\|^{2}+\varkappa(x)x^{2}+x\left(m_{1}(z)\|z\|\right.\\
\left.+m_{2}(x)|x|-\rho(x)x+\varpi\right)\leq-\|\xi\|^{2}+\varpi^{2}
\end{gather*}

Now, let us examine the $\varpi$-dynamics. Note that $\dot{\varpi}=\frac{d\kappa(x(t),z(t))}{dt}=-\pi(x)\dot{x}$
where $\tau(x)=\frac{d\rho(x)}{dx}x+\rho(x)$. Let $\bar{\tau}(x)=|\tau(x)|$
and decompose $\bar{\tau}(x)$ to be $\bar{\tau}(x)=\tau_{x}(x)+\tau_{c}$
where $\tau_{x}(0)=0$ and $\tau_{c}\geq0$. As a result, $\tau_{x}(x)\leq m_{3}(x)|x|$
for some   function $m_{3}\in\mathcal{SN}$. Let $U(\varpi)=\frac{1}{2}\varpi^{2}$.
We claim that
\begin{equation}
\dot{U}(\varpi)\leq\alpha_{\varpi}(\|\varpi\|)+\sigma_{\varpi}(\|\xi\|)\label{eq:Uw}
\end{equation}
with some functions $\alpha_{\varpi},\sigma_\varpi\in\mathcal{K}_\infty$.
In fact, one possible calculation of the derivative of $U(\varpi)$
is given below
\begin{gather*}
\dot{U}(\varpi)\leq|\varpi|\left[|\tau(x)|(m_{1}(z)\|z\|+m_{2}(x)|x|+\rho(x)|x|+|\varpi|)\right]\\
\leq|\varpi|\left(\tau_{c}m_{1}(z)\|z\|+\tau_{x}(x)m_{1}(z)\|z\|+\bar{\tau}(x)m_{2}(x)|x|\right.\\
\left.+\bar{\tau}(x)\rho(x)|x|+\tau_{c}|\varpi|+\tau_{x}(x)|\varpi|\right)\\
\leq\tau_{c}^{2}m_{1}^{2}(z)\|z\|^{2}/4+m_{1}^{4}(z)\|z\|^{4}/8+m_{3}^{4}(x)x{}^{4}/8\\
+\bar{\tau}^{2}(x)(m_{2}(x)+\rho(x))^{2}x{}^{2}/4+m_{3}^{2}(x)x{}^{2}/4\\
+(3+\tau_{c})\varpi^{2}+\varpi^{4}\\
\leq\bar{\alpha}_{\varpi}(\varpi)\|\varpi\|^{2}+\bar{\sigma}(\xi)\|\xi\|^{2}
\end{gather*}
where $\bar{\sigma}(\xi)$ and $\bar{\alpha}_{\varpi}(\varpi)$ are
selected as 
\begin{gather*}
\bar{\alpha}_{\varpi}(\varpi)\geq(3+\tau_{c})+\varpi^{2},\\
\bar{\sigma}(\xi)\geq\tau_{c}^{2}m_{1}^{2}(z)/4+m_{1}^{4}(z)\|z\|^{2}/8+m_{3}^{4}(x)x^{2}/8\\
+\bar{\tau}^{2}(x)(m_{2}(x)+\rho(x))^{2}/4+m_{3}^{2}(x)/4.
\end{gather*}
Let us choose $\alpha_{\varpi}(\|\varpi\|)\geq\bar\alpha_{\varpi}(\varpi)\|\varpi\|^2$ and $\sigma_{\varpi}(\|\xi\|)\geq\bar\sigma(\xi)\|\xi\|^2$. Then, (\ref{eq:Uw}) is satisfied.

Let function $\alpha_{q}(\|\xi\|)\geq2\sigma_{\varpi}(\|\xi\|)+\|\xi\|^{2}$.
Since there exists
an ISS Lyapunov function $V_{\xi}(\xi)$ for (\ref{eq:rb_sys_cl}),
by changing supply function method, there exists another ISS Lyapunov
function $V_{q}(\xi)\sim\{\underline{\alpha}_{q},\bar{\alpha}_{q},\alpha_{q},\sigma_{q}\mid\eqref{eq:rb_sys_cl}\}$
for some $\mathcal{K}_{\infty}$ functions $\underline{\alpha}_{q}$,
$\bar{\alpha}_{q}$, and $\sigma_{q}$, that are calculated accordingly.
Moreover, $\limsup_{s\rightarrow0^{+}}\sigma_q(s)/s^{2}>0$.  Select $\mathcal{K}_\infty$ function $\gamma$ as
\begin{equation}
\gamma(s)\geq\sigma_{q}(s)/2+\alpha_{\varpi}(s)+\frac{1}{2}s^2 .\label{eq:gamma}
\end{equation}
Let $\bar V(\xi,\varpi)=V_q(\xi)+U(\varpi))$ be the  Lyapunov function candidate for the system  (\ref{eq:rb_sys_cl}). We can check that  conditions of Theorem \ref{lem:uniform_lem} is satisfied.
Moreover, 
\begin{gather}
\bar{V}(x(t),\varpi(t))-\bar{V}(x(t_{k}^+),\varpi(t_{k}^+))\nonumber \\
\leq-\int_{t_{k}^+}^{t}\|\xi(\tau)\|^2d\tau-\int_{t_{k}^+}^{t}\varpi^2(\tau) d\tau, \forall t\in[t_{k}^+,t_{k+1}]
\end{gather}
Applying Theorem \ref{lem:uniform_lem}  completes the proof.
\eproof 

As has been done in Section \ref{sec:SCD}, the periodic sampling
law can be found by following the idea of Theorem \ref{thm:robust}.
The proof is straightforward and omitted here.


\bproposition \label{thm:ps_robust} Consider the system composed of (\ref{eq:rb_sys})
and (\ref{eq:rb_con}) under Assumption \ref{ass:z}. 
Let functions $V_q$, $\gamma$, $\alpha_{\varpi}$, $\sigma_{\varpi}$  and $\underline{\alpha}_{q}$ 
be defined in the proof of Theorem \ref{thm:robust}.
Let $\mathcal{X}:=\{\xi\in\mathbb{R}^{p+1}\mid\|z\|\leq z_c,\; x=x_0\}$ and
$R_{0}=\max_{\xi\in{\mathcal{X}}}\{ \alpha_{\varpi}^{-1}\circ\sigma_{\varpi}\circ\underline{\alpha}_{q}{}^{-1}\left(V_{q}(\xi)\right)\}$. 
Find $T$ such that 
\begin{gather}
2T\gamma(s)<\underline{\alpha}_{\varpi}(s),\;\forall0<s\leq R_{0},\label{eq:robust_t}
\end{gather}
the equilibrium point $\xi=0$ of the system is asymptotically stable asymptotically stable for any initial condition $\xi(t_{0})\in \mathcal{X}$.
\eproposition

The discussion in Remark \ref{rem:psc_init} shows that the initial
condition of the signal $\xi$ must be known for the design of periodic
sampling. Although the signal $z$ is not available for the feedback,
Proposition \ref{thm:ps_robust} requires the initial condition of the
signal $z$ or at least its bound be known, which is used to estimate
$T$ in (\ref{eq:robust_t}).  

\brem In \cite{Liu2015,Liu2015a}, the decay rate of immeasurable states $z$ is estimated by an auxiliary dynamic system  and then used  for constructing the event-triggered law.  In comparison, a static
sampling controller is proposed in Theorem \ref{thm:robust},
 that  is easier to design and implement in practice.
Moreover, our method can also be used to derive the periodic sampling
control law when the initial bound of the signal $z$ is known. \erem

\section{Even-triggered Adaptive Control \label{sec:adaptive}}

In this section, we will solve two types of classical adaptive control
problem in the event-triggered setting exploiting the sampling control
scheme proposed in Section \ref{sec:SCD}.

\subsection{Adaptive Control with Uncertain Parameters \label{sub:acup}}

We consider the event-triggered adaptive stabilization problem of
a class of nonlinear systems with unknown parameters, as follows,
\begin{eqnarray}
\dot{x} & = & f\t(x)\theta+u\label{eq:adaptive_sys}
\end{eqnarray}
where $x\in\mathbb{R}$ is the state and $\theta\in\mathbb{R}^{l}$
is an unknown constant parameter vector. Note that function $f(x):\mathbb{R}\mapsto \mathbb{R}^{l}$ does
not necessarily vanish at $x=0$. Without loss of generality, we assume
$f(0)\neq0$ and $\theta$ is bounded.

\bass \label{ass:theta_b} (1) The unknown parameter $\theta$ is
bounded with a known bound $\theta_{b}$, i.e., $\|\theta\|\leq\theta_{b}$;
(2) $f(0)\neq0$. \eass

In this case, the argument in Remark \ref{rem:bcon} does not apply
and inequality (\ref{eq:bcon}) in Assumption \ref{ass:SP-Bound}
might not hold. Therefore, the development of the controller design
in Section \ref{sec:SCD} must be modified to suit the problem. Similar
to the continuous-time adaptive control for the system (\ref{eq:adaptive_sys}),
we propose the event-triggered controller as follows,
\begin{eqnarray}
u & = & \kappa(x(t_{k}),\hat{\theta}(t_{k})),\;t\in[t_{k},t_{k+1}),\;k\in\mathbb{N}\nonumber \\
\dot{\hat{\theta}} & = & \Lambda x\rho(x)+\Lambda\varsigma(\hat{\theta},\varpi)\label{eq:ada_con}
\end{eqnarray}
where $\kappa(x,\hat{\theta})=-f\t(x)\hat{\theta}-5x/4$, $\Lambda>0\in\mathbb{R}^{l\times l}$  is a diagonal matrix and $\rho(x)$, functions
$\varsigma(\hat{\theta},\varpi)$ are to be designed. Note that the
term $\varsigma(\hat{\theta},\varpi)$ does not appear in traditional
continuous-time adaptive control and is introduced particularly for
the event-triggered control. In this paper, we do not consider the
challenging case where the adaptation dynamics $\hat{\theta}$ is
sampled and it remains our future research direction. Define the sampling
error $\varpi(t)$ as $\varpi(t)=\kappa(x(t_{k}),\hat{\theta}(t_{k}))-\kappa(x(t),\hat{\theta}(t)),\;t\in[t_{k},t_{k+1}).$
Then, the closed-loop system is rewritten as 
\begin{eqnarray}
\dot{x} & = & -f\t(x)\tilde{\theta}-\frac{5}{4}x+\varpi,\nonumber \\
\dot{\hat{\theta}} & = & \Lambda x\rho(x)+\Lambda\varsigma(\hat{\theta},\varpi),\nonumber \\
\dot{\varpi} &=& \frac{\partial \kappa(x,\hat{\theta})}{\partial x} \dot{x}+\frac{\partial \kappa(x,\hat{\theta})}{\partial \hat{\theta}}\dot{\hat{\theta}},
\label{eq:adaptive_cl}
\end{eqnarray}
where $\tilde{\theta}=\hat{\theta}-\theta$ is estimation error. Based on the event-triggered
controller design in Section \ref{sec:SCD}, the following theorem
is obtained.

\bthm \label{thm:adap} Suppose the system composed of (\ref{eq:adaptive_sys})
and (\ref{eq:ada_con}) satisfies Assumption \ref{ass:theta_b}. Then,
there exists smooth positive functions $\rho(x):\mathbb{R}\mapsto\mathbb{R}^{l}$,
$\varsigma(\hat{\theta},\varpi):\mathbb{R}^{l}\times\mathbb{R}\mapsto\mathbb{R}$
for the controller (\ref{eq:ada_con}) and $\bar{\gamma}(x,\hat{\theta},\varpi):\mathbb{R}\times\mathbb{R}^{l}\times\mathbb{R}\mapsto\mathbb{R}_{>0}$
such that the event-triggered law is designed
as
\begin{gather}
t_{k+1}=\inf_{t\geq t_{k}}\{2(t-t_{k})\max_{\tau\in[t_{k},t]}\{\bar{\gamma}(x(\tau),\hat{\theta}(\tau),\varpi(\tau))\varpi^{2}(\tau)\},\nonumber \\
\geq\frac{1}{2}\|\varpi_{[t_{k},t]}\|^{2},\;\text{and }\|\varpi_{[t_{k},t]}\|\neq0\}\label{eq:sd_con_ada}
\end{gather}
then Zeno-behavior is avoided and the equilibrium point $x=0$ is
globally asymptotically stable. Moreover, $\lim_{t\rightarrow\infty}\hat{\theta}(t)=\theta.$
The algorithm of event-triggered adaptive controller design is summarized
in Algorithm \ref{alg:sdac}. \ethm \proofnow First, let us consider
the $\varpi$-dynamics. Denote $\tau(x,\hat{\theta})=\frac{\partial f\t(x)}{\partial x}\hat{\theta}+\frac{5}{4}$
and $U(\varpi)=\frac{1}{2}\varpi^{2}.$ Note that
\begin{gather*}
\dot{\varpi}=\frac{d\kappa(x(t),\theta(t))}{dt}=-\left(\tau(x,\hat{\theta})\dot{x}+f\t(x)\dot{\hat{\theta}}\right)\\
=\tau(x,\hat{\theta})f\t(x)\tilde{\theta}+\tau(x,\hat{\theta})\frac{5}{4}x-\tau(x,\hat{\theta})\varpi\\
-f\t(x)\Lambda x\rho(x)-f\t(x)\Lambda\varsigma(\hat{\theta},\varpi).
\end{gather*}
Decompose $f(x)$ as $f(x)=\bar{\varrho}_{x}(x)+\bar{\varrho}_{c}$
where $\bar{\varrho}_{x}(x)$ depends on $x$ satisfying $\bar{\varrho}_{x}(0)=0$
and $\bar{\varrho}_{c}$ is a constant vector. Decompose $\tau(x,\hat{\theta})$
as $\tau(x,\hat{\theta})=\bar{\tau}_{x}\t(x)\hat{\theta}+\bar{\tau}_{c}(\hat{\theta})$
where $\bar{\tau}_{x}(x)$ depends on $x$ satisfying $\bar{\tau}_{x}(0)=0$
and $\bar{\tau}_{c}(\hat{\theta})$ is a scalar  possibly depending on $\hat{\theta}$.
Denote $\varrho_{x}(x)=\|\bar\varrho_{x}(x)\|$, $\varrho_{c}=\|\bar{\varrho}_{c}\|$,
$\tau_{x}(x)=|\bar{\tau}_{x}\t(x)|$, $\tau_{c}(\hat{\theta})=|\bar{\tau}_{c}(\hat{\theta})|$
and $\chi(x)=\|f\t(x)\Lambda\|$. Note that $\|\tilde{\theta}\|^{2}\leq2\|\theta\|^{2}+2\|\hat{\theta}\|^{2}\leq2\theta_{b}^{2}+2\|\hat{\theta}\|^{2}$.
As a result,
\begin{gather}
\dot{U}(\varpi)=\tau(x,\hat{\theta})f\t(x)\tilde{\theta}\varpi+\frac{5}{4}\tau(x,\hat{\theta})x\varpi\nonumber \\
-\tau(x,\hat{\theta})\varpi^{2}-f\t(x)\Lambda\rho(x)x\varpi-f\t(x)\Lambda\varsigma(\hat{\theta},\varpi)\varpi\nonumber \\
\leq\left[\tau_{x}(x)\|f(x)\|\|\hat{\theta}\|+\varrho_{x}(x)\tau_{c}(\hat{\theta})\right]\|\tilde{\theta}\||\varpi|+\bar{\varrho}_{c}\t\bar{\tau}_{c}(\hat{\theta})\tilde{\theta}\varpi\nonumber \\
+\left(\frac{5}{4}\tau_{x}(x)\|\hat{\theta}\|+\frac{5}{4}\tau_{c}(\hat{\theta})\right)|x||\varpi|-\tau(x,\hat{\theta})\varpi^{2}\nonumber \\
+\chi^{2}(x)x^{2}/4+|\rho(x)\varpi|^{2}-f\t(x)\Lambda\varsigma(\hat{\theta},\varpi)\varpi\nonumber \\
\leq\tau_{x}^{2}(x)\|f(x)\|^{2}/4+\varrho_{x}^{2}(x)/4+\frac{25}{64}\tau_{x}^{2}(x)x^{2}\nonumber \\
+\frac{25}{64}x^{2}+\chi^{2}(x)x^{2}/4+\|\hat{\theta}\|^{2}\|\tilde{\theta}\|^{2}\varpi^{2}\nonumber \\
+\tau_{c}^{2}(\hat{\theta})\|\tilde{\theta}\|^{2}\varpi^{2}+\|\hat{\theta}\|^{2}\varpi^{2}-\tau(x,\hat{\theta})\varpi^{2}+|\rho(x)|^{2}\varpi{}^{2}\nonumber \\
+\bar{\varrho}_{c}\t\bar{\tau}_{c}(\hat{\theta})\tilde{\theta}\varpi-f\t(x)\Lambda\varsigma(\hat{\theta},\varpi)\varpi\label{eq:Uw_dot}
\end{gather}
Since $\tau_{x}(0)=0$ and $\varrho_{x}(0)=0$, $\tau_{x}(x)\leq m_{1}(x)|x|$
and $\varrho_{x}(x)\leq m_{2}(x)|x|$ for some   functions
$m_{1},m_{2}\in\mathcal{SN}$. Select two smooth positive functions $\bar{\sigma}_{\varpi}$
and $\bar{\alpha}_{\varpi}$ as 
\begin{gather}
\bar{\sigma}_{\varpi}(x)\geq\frac{1}{4}m_{1}^{2}(x)\|f(x)\|^{2}+\frac{1}{4}m_{2}^{2}(x)+\frac{25}{64}m_{1}^{2}(x)x^{2}\nonumber \\
+\frac{25}{64}+\chi^{2}(x)/4\nonumber \\
\bar{\alpha}_{\varpi}(x,\hat{\theta})\geq\left(\|\hat{\theta}\|^{2}+\tau_{c}^{2}(\hat{\theta})\right)\left(2\theta_{b}^{2}+2\|\hat{\theta}\|^{2}\right)\nonumber \\
+\|\hat{\theta}\|^{2}+|\rho(x)|^{2}-\tau(x,\hat{\theta})\label{eq:sig_alpha}
\end{gather}
Moreover, $\bar{\alpha}_{\varpi}(x,\hat{\theta})$ is bounded when
signals $x$ and $\hat{\theta}$ are bounded. As a result, it follows
from (\ref{eq:Uw_dot}) that 
\begin{gather}
\dot{U}(\varpi)\leq\bar{\sigma}_{\varpi}(x)x{}^{2}+\bar{\alpha}_{\varpi}(x,\hat{\theta})\varpi{}^{2}\nonumber \\
+\bar{\varrho}_{c}\t\bar{\tau}_{c}(\hat{\theta})\tilde{\theta}\varpi-f\t(x)\Lambda\varsigma(\hat{\theta},\varpi)\varpi.\label{eq:Uw_dot2}
\end{gather}
It should be noted that the selection of $\bar{\sigma}_{\varpi}(x)$
only depends on the function $f(x)$ and not on function $\rho(x)$
in $\hat{\theta}$-dynamics.  In other words, $\rho(x)$ can be selected
to be dependent on $\bar{\sigma}_{\varpi}(x)$.

Then, let us consider the $x$-dynamics. Let $V_{x}(x)=\frac{1}{2}x^{2}$.
The derivative of $V_{x}(x)$ along the $x$-dynamics is 
\[
\dot{V}_{x}(x)=x\left(-f\t(x)\tilde{\theta}-\frac{5}{4}x+\varpi\right)\leq-x^{2}+\varpi^{2}-xf\t(x)\tilde{\theta}.
\]
Following changing supply function lemma (in Lemma 2.5 of \cite{Chen2015Book}),
for $\Delta(x)=2\bar{\sigma}_{\varpi}(x)+1$, there exists a non-decreasing
function $\rho_{q}$ such that 
\begin{equation}
\frac{1}{2}\rho_{q}\left(\frac{1}{2}x^2\right)\geq\Delta(x)\label{eq:rho_q}
\end{equation}
and a new supply function $V_{x}'(x)=\int_{0}^{V(x)}\rho_{q}(s)ds$
such that its derivative along the $x$-dynamics satisfies
\begin{gather}
\dot{V}_{x}'(x)\leq-\Delta(x)x^{2}+\sigma(\varpi)\|\varpi\|^{2}-\rho_{q}(V(x))xf\t(x)\tilde{\theta}\nonumber \\
\leq-\Delta(x)x{}^{2}+\sigma(\varpi)\|\varpi\|^{2}-\rho_{q}\left(\frac{1}{2}x^{2}\right)xf\t(x)\tilde{\theta}\label{eq:Vx'}
\end{gather}
for some   $\sigma\in\mathcal{SP}$. Since $\bar{\sigma}_{\varpi}(x)$
is not dependent on $\rho(x)$, so is function $\rho_{q}(x)$. Then,
one can choose $\rho(x)$ and $\varsigma(\hat{\theta},\varpi)$ in
(\ref{eq:ada_con}) to be 
\begin{gather}
\rho(x)=\rho_{q}\left(\frac{1}{2}x^{2}\right)f(x)\nonumber \\
\varsigma(\hat{\theta},\varpi)=-\bar{\varrho}_{c}\bar{\tau}_{c}(\hat{\theta})\varpi.\label{eq:sel}
\end{gather}

Now, let us consider the $\tilde{\theta}$-dynamics. Let $V_{\theta}(\tilde{\theta})=\tilde{\theta}\t\Lambda^{-1}\tilde{\theta}/2$
and its derivative along $\tilde{\theta}$-dynamics becomes
\[
\dot{V}_{\theta}(\tilde{\theta})=\rho_{q}\left(\frac{1}{2}x^{2}\right)xf\t(x)\tilde{\theta}-\tilde{\theta}\t\bar{\varrho}_{c}\bar{\tau}_{c}(\hat{\theta})\varpi.
\]
Denote $\xi=\text{col}(x,\tilde{\theta})$. Letting $V(\xi)=V_{x}'(x)+V_{\theta}(\tilde{\theta})$
leads to 
\begin{equation}
\dot{V}(\xi)=-\Delta(x)x{}^{2}+\sigma(\varpi)\|\varpi\|^{2}-\tilde{\theta}\t\bar{\varrho}_{c}\bar{\tau}_{c}(\hat{\theta})\varpi.
\end{equation}
By the selection of $\rho(x)$ and $\varsigma(\hat{\theta},\varpi)$
in (\ref{eq:sel}), the inequality (\ref{eq:Uw_dot2}) further becomes
\begin{equation}
\dot{U}(\varpi)\leq\bar{\sigma}_{\varpi}(x)x^{2}+\tilde{\alpha}_{\varpi}(x,\hat{\theta})\varpi{}^{2}+\bar{\varrho}_{c}\t\bar{\tau}_{c}(\hat{\theta})\tilde{\theta}\varpi\label{eq:Uw_dot_ada}
\end{equation}
where 
\begin{equation}
\tilde{\alpha}_{\varpi}(x,\hat{\theta})\geq\bar{\alpha}_{\varpi}(x,\hat{\theta})+\|f\t(x)\Lambda\bar{\varrho}_{c}\bar{\tau}_{c}(\hat{\theta})\|\label{eq:alpha_w}
\end{equation}

Finally, the event-triggered law design is presented. Let 
\begin{equation}
\bar{\gamma}(x,\hat{\theta},\varpi)=\sigma(\varpi)/2+\tilde{\alpha}_{\varpi}(x,\hat{\theta})+\frac{1}{2}\label{eq:gamma_f}
\end{equation}
where $\sigma(\varpi)$ is given in (\ref{eq:Vx'}).
One has 
\begin{gather*}
U(\varpi(t))=\int_{t_{k}^+}^{t}\dot{U}(\varpi(\tau))d\tau\\
\leq\int_{t_{k}^+}^{t}\left(\bar{\sigma}_{\varpi}(x(\tau))x^{2}(\tau)+\bar{\varrho}_{c}\t\bar{\tau}_{c}(\hat{\theta}(\tau))\tilde{\theta}(\tau)\varpi(\tau)\right)d\tau\\
+\int_{t_{k}^+}^{t}\tilde{\alpha}_{\varpi}(x(\tau),\hat{\theta}(\tau))\varpi{}^{2}(\tau)d\tau,\;t\in[t_{k},t_{k+1}]
\end{gather*}
which together with event-triggered law (\ref{eq:sd_con_ada}) implies
\begin{gather*}
U(\varpi(t))\leq2\int_{t_{k}^+}^{t}\left(\bar{\sigma}_{\varpi}(x(\tau))x^{2}(\tau)+\bar{\varrho}_{c}\t\bar{\tau}_{c}(\hat{\theta}(\tau))\tilde{\theta}(\tau)\varpi(\tau)\right.\\
\left.-\frac{1}{2}\sigma(\varpi(\tau))\varpi{}^{2}(\tau)-\frac{1}{2}\varpi^{2}(\tau)\right)d\tau,\;t\in[t_{k},t_{k+1}]
\end{gather*}
Let $\bar{V}(x,\tilde{\theta},\varpi)=V_{x}'(x)+V_{\theta}(\tilde{\theta})+U(\varpi)$
be the Lyapunov function candidate for the closed-loop system (\ref{eq:adaptive_cl}).
As a result,
\begin{gather*}
\bar{V}(x(t),\varpi(t))-\bar{V}(x(t_{k}^+),\varpi(t_{k}^+))\leq\\
\leq-\int_{t_{k}^+}^{t}x^{2}(\tau)d\tau-\int_{t_{k}^+}^{t}\varpi^{2}(\tau)d\tau,\;\forall t\in[t_{k},t_{k+1}]
\end{gather*}
Therefore, $\bar{V}(x,\tilde{\theta},\varpi)$ is monotonically decreasing
and all signal $x$, $\hat{\theta}$ and $\varpi$ are bounded. So,
there exists a positive constant $\gamma_{c}$ depending on the bound
of $x$, $\hat{\theta}$ and $\varpi$ such that $\bar{\gamma}(x,\hat{\theta},\varpi)\leq\gamma_{c}$.
Moreover, there exists a $t_{c}>0$ such that 
\begin{gather*}
2(t-t_{k})\max_{\tau\in[t_{k},t]}\{\bar{\gamma}(x(\tau),\hat{\theta}(\tau),\varpi(\tau))\varpi^{2}(\tau)\},\\
\leq2(t-t_{k})\gamma_{c}\|\varpi_{[t_{k},t]}\|^{2}\\
\leq\frac{1}{2}\|\varpi_{[t_{k},t]}\|^{2},\;t\leq t_{k}+t_{c}
\end{gather*}
and noting (\ref{eq:sd_con_ada}) shows $t_{k+1}\geq t_{k}+t_{c}.$
Therefore, it is Zeno free. 

The convergence of the signal $\text{col}(x(t),\varpi(t))$ to zero
follows similar argument in the proof of Theorem \ref{lem:uniform_lem}.
When $x\equiv0$ and $\varpi\equiv0$, it shows that $\lim_{t\rightarrow\infty}f\t(x(t))\tilde{\theta}(t)=0$.
Because $f(0)$ is PE, it is proved that $\lim_{t\rightarrow\infty}\tilde{\theta}(t)=0$
or $\lim_{t\rightarrow\infty}\hat{\theta}(t)=\theta$ by Lemma 2.4
of \cite{Chen2015Book}. \eproof 

\begin{algorithm}
1. Select smooth positive functions $\bar{\sigma}_{\varpi}$ and $\bar{\alpha}_{\varpi}$
to satisfy (\ref{eq:sig_alpha}); 

2. Find a positive function $\rho_{q}$ to satisfy (\ref{eq:rho_q})
where $\Delta(x)=2\bar{\sigma}_{\varpi}(x)+1$ and $\sigma(\varpi)$
to satisfy (\ref{eq:Vx'});

3. Choose $\rho(x)$ and $\varsigma(\hat{\theta},\varpi)$ to be in
(\ref{eq:sel});

4. Find a positive function $\tilde{\alpha}_{\varpi}(x,\hat{\theta})$
to satisfy (\ref{eq:alpha_w}) and $\bar\gamma(x,\hat{\theta},\varpi)$
to be in (\ref{eq:gamma_f});

5. Construct the event-triggered law in (\ref{eq:sd_con_ada}).

\caption{Event-triggered Adaptive Control\label{alg:sdac}}
\end{algorithm}

\brem In \cite{huang2019}, the event-triggered adaptive control
problem is solved when the global Lipschitz condition is assumed on
the function $f(x)$, while it is not required in our method. It would
be very interesting to consider that the adaptation dynamics is also
sampled as in \cite{huang2019} in our future work. Compared with the
development in Section \ref{sec:SCD}, the derivative of $U(w)$ is
upper bounded as in (\ref{eq:Uw_dot2}) which is more complicated
than  (\ref{eq:bcon}) in Assumption \ref{ass:SP-Bound}. Therefore,
the event-triggered law in (\ref{eq:sd_con_ada}) is also more involved and
not merely determined by $\varpi$ but also $x$ and $\hat{\theta}$.
\erem

\subsection{Event-triggered Stabilization with Dynamic Gain}

In this subsection, we consider the dynamic system 
\begin{eqnarray}
\dot{z} & = & q(z,x,d)\nonumber \\
\dot{x} & = & f(z,x,d)+bu,\;b>0\label{eq:dg_sys}
\end{eqnarray}
where $z\in\mathbb{R}^{p}$ and $x\in\mathbb{R}$ are state variables,
$u\in\mathbb{R}$ is the input and $d\in\mathbb{D}$ is the uncertainties
belonging to a compact set $\mathbb{D}\in\mathbb{R}^{l}$. The functions
$q$ and $f$ are sufficiently smooth with $q(0,0,d)=0$ and $f(0,0,d)=0$
for all $d\in\mathbb{D}$. In contrast to Section \ref{sec:robust},
the size of $\mathbb{D}$ and controller gain $b$ are not known and a dynamic gain technique
is required. The state $z$ is also assumed to be not available for
feedback control and thus $z$-dynamics is  dynamic uncertainty.
A few Assumptions are list as follows.

\bass \label{ass:b_bound} $b\leq\bar{b}$ for some known constant
$\bar{b}$. \eass 

\bass \label{ass:za}The $z$-subsystem in (\ref{eq:rb_sys}) has
an ISS-Lyapunov function $V(z)\sim\{\underline{\alpha},\bar{\alpha},\alpha,\hat{\sigma}\mid\dot{z}=q(z,x,d)\}$
and functions $\alpha(s)$ and $\sigma(s)$ satisfy $\limsup_{s\rightarrow0^{+}}s^{2}/\alpha(s)<\infty,\;\limsup_{s\rightarrow0^{+}}\hat{\sigma}(s)/s^{2}<\infty.$
Moreover, the functions $\underline{\alpha},\bar{\alpha},\alpha$
are known and the function $\hat{\sigma}$ is known up to a constant
factor in the sense that there exist an unknown constant $p$ and
a known function \textgreek{sv} such that $\hat{\sigma}=p\sigma$.
\eass

Note that the system (\ref{eq:dg_sys}) differs from that studied
in Section \ref{sec:robust} in some aspects. First, in Section \ref{sec:robust},
since the controller gain is unity, we can apply sufficient high gain
specified by $\rho(x)$ to dominate uncertainties and stabilize the
system. Here, the controller gain $b$ is not known. Although the
upper bound of $b$ is known, it is still not possible to calculate
the controller gain that is sufficiently high to stabilize the system.
Second, although the $z$-subsystem in (\ref{eq:dg_sys}) admits an
ISS-Lyapunov function, the input gain $\hat{\sigma}$ is only known
to a constant factor. Third, since the size of $\mathbb{D}$ is unknown,
by Corollary 11.1 in \cite{Chen2015Book}, there exists a positive
number $c$, which depends on the size of $\mathbb{D}$ and is also
unknown, and two positive and sufficiently smooth known functions
$m_{1}$ and $m_{2}$, such that 
\begin{equation}
|f(z,x,d)|\leq cm_{1}(z)\|z\|+cm_{2}(x)|x|,\;\forall d\in\mathbb{D}.\label{eq:f_bound_2}
\end{equation}
These three differences call for the dynamic gain stabilization technique.

For the continuous-time stabilization of the system (\ref{eq:dg_sys}),
dynamic gain controller using universal adaptive control technique
is proposed in \cite{Chen2015Book}. It takes the form of 
\begin{eqnarray}
u & = & \kappa(x,\theta)=-\hat{\theta}\rho(x)x\nonumber \\
\dot{\hat{\theta}} & = & \lambda\rho(x)x^{2},\;\lambda>0\label{eq:dg_con}
\end{eqnarray}
where $\rho(x)$ is to be specified. The second equation of (\ref{eq:dg_con})
is the adaptation dynamics for dynamic gain. For the event-triggered control,
we adopt the method developed in Section \ref{sec:SCD} and propose
the controller as follows 
\begin{align}
u(t) & =\kappa(x(t_{k}),\hat{\theta}(t_{k})),\;t\in[t_{k},t_{k+1}),\;k\in\mathbb{N} \nonumber \\
\dot{\hat{\theta}} & =  \lambda\rho(x)x^{2},\;\lambda>0 \label{eq:dg}
\end{align}
Define the sampling error $\varpi(t)$ as $\varpi(t)=\kappa(x(t_{k}),\theta(t_{k}))-\kappa(x(t),\theta(t)),\;t\in[t_{k},t_{k+1}).$
Then, the event-triggered closed-loop system is rewritten as 
\begin{eqnarray}
\dot{z} & = & q(z,x,d)\nonumber \\
\dot{x} & = & f(z,x,d)-b\hat{\theta}\rho(x)x+b\varpi,\nonumber \\
\dot{\hat{\theta}} & = & \lambda\rho(x)x^{2},\;\lambda>0\nonumber·\\
\dot{\varpi} &=& \frac{\partial \kappa(x,\hat{\theta})}{\partial x} \dot{x}+\frac{\partial \kappa(x,\hat{\theta})}{\partial \hat{\theta}}\dot{\hat{\theta}},
\label{eq:gain_cl}
\end{eqnarray}

\bthm \label{thm:dg} Suppose the system composed of (\ref{eq:dg_sys}) and (\ref{eq:dg})
satisfies Assumption \ref{ass:b_bound} and \ref{ass:za}. Then, there
exist smooth positive functions $\rho(x):\mathbb{R}\rightarrow\mathbb{R}_{>0}$,
and $\bar{\gamma}(x,\hat{\theta},\varpi):\mathbb{R}\times\mathbb{R}\times\mathbb{R}\rightarrow\mathbb{R}_{>0}$
such that event-triggered law ) is designed
as
\begin{gather}
t_{k+1}=\inf_{t\geq t_{k}}\{2(t-t_{k})\max_{\tau\in[t_{k},t]}\{\bar{\gamma}(x(\tau),\hat{\theta}(\tau),\varpi(\tau))\varpi^{2}(\tau)\}\nonumber \\
\geq\frac{1}{2}\|\varpi_{[t_{k},t]}\|^{2},\;\text{and }\|\varpi_{[t_{k},t]}\|\neq0\},\label{eq:sd_con_ada-1}
\end{gather}
then the equilibrium
point $\text{col}(z,x)=0$ of the system is globally asymptotically
stable and Zeno-behavior is avoided.  \ethm \proofnow Denote $\xi=\text{col}(z,x)$ and $\Xi=\text{col}(z,x,\hat{\theta})$.\textcolor{red}{{}
}First, let us examine $\Xi$-dynamics. By changing supply function
method, for $\Delta(z)=m_{1}^{2}(z)+1$, there exists another ISS
Lyapunov function, $V_{z}(z)\sim\{\underline{\alpha}_{z},\bar{\alpha}_{z},\Delta(z)\|z\|^{2},\varkappa(x)x^{2}\mid\dot{z}=q(z,x,d)\}$
for some $\mathcal{K}_{\infty}$ functions $\underline{\alpha}_{z}$,
$\bar{\alpha}_{z}$, and a smooth function $\varkappa$, that are
calculated accordingly. Let $\tilde{\theta}=\hat{\theta}-\theta$
where $\theta$ is a positive number to be specified later. Define
the Lyapunov function candidate as 
\begin{equation}
V(z,x,\hat{\theta})=V_{z}(z)+x^{2}/2+b\tilde{\theta}^{2}/(2\lambda)\label{eq:Vzxt}
\end{equation}
The calculation of the derivative of $V(z,x,\tilde{\theta})$, along
the trajectory of (\ref{eq:gain_cl}), shows that 
\begin{gather}
\dot{V}(z,x,\tilde{\theta})\leq-\Delta(z)\|z\|^{2}+\varkappa(x)x^{2}+x(cm_{1}(z)\|z\|\nonumber \\
+cm_{2}(x)|x|-b\theta\rho(x)x+b\varpi)\nonumber \\
\leq-\|z\|^{2}+(\varkappa(x)+c^{2}/4+cm_{2}(x)+b^{2}q/4\nonumber \\
-b\theta\rho(x))x^{2}+\varpi^{2}/q\label{eq:V_dot_LTS}
\end{gather}
where $q$ is to be specified. Let $\rho(x)$ be

\[
\rho(x)\geq\max\{1,\varkappa(x),m_{2}(x)\}
\]
and 
\[
\theta\geq(1+c^{2}/4+c)/b+bq/4.
\]
As a result, one has 
\begin{equation}
\dot{V}(z,x,\hat{\theta})\leq-\|\xi\|^{2}+\varpi^{2}/q.\label{eq:V_dot}
\end{equation}

Next, the $\varpi$-dynamics is examined as follows
\begin{gather*}
\dot{\varpi}=\frac{d\kappa(x(t),z(t))}{dt}=-\dot{\hat{\theta}}\rho(x)x-\hat{\theta}\frac{d\rho(x)}{dx}x\dot{x}-\hat{\theta}\rho(x)\dot{x}\\
=-\lambda\rho^{2}(x)x^{3}-\hat{\theta}\pi(x)\left(f(\psi,x,d)-b\hat{\theta}\rho(x)x+b\varpi\right)
\end{gather*}
where $\pi(x)=\frac{d\rho(x)}{dx}x+\rho(x)$. Let $\bar{\tau}(x)=|\tau(x)|$
and decompose $\bar{\tau}(x)$ to be $\bar{\tau}(x)=\tau_{x}(x)+\tau_{c}$
where $\tau_{x}(0)=0$ and $\tau_{c}\geq0$. As a result, $\tau_{x}(x)\leq m_{3}(x)|x|$
for some  function $m_{3}\in\mathcal{SN}$. Let $U(\varpi)=\frac{1}{2}\varpi^{2}$
and one has
\begin{gather*}
\dot{U}(\varpi)\leq|\varpi||\dot{\varpi}|\leq|\varpi|\left(\lambda\rho^{2}(x)|x|^{3}+b\hat{\theta}^{2}\bar{\pi}(x)\rho(x)|x|\right.\\
+b\bar{\tau}(x)|\hat{\theta}||\varpi|+c\hat{\theta}\bar{\pi}(x)m_{2}(x)|x|\\
\left.+c^{2}(\hat{\theta}\tau_{x}(x))^{2}/4+c\hat{\theta}\tau_{c}m_{1}(z)\|z\|+m_{1}^{2}(z)\|z\|^{2}\right)\\
\leq\lambda^{2}\rho^{4}(x)x{}^{6}/4+b^{2}\bar{\pi}^{2}(x)\rho^{2}(x)x{}^{2}/4+c^{2}\bar{\pi}^{2}(x)m_{2}^{2}(x)x{}^{2}/4\\
+c^{4}\tau_{x}^{4}(x)/64+m_{1}^{4}(z)\|z\|^{4}/4+c^{2}\tau_{c}^{2}m_{1}^{2}(z)\|z\|^{2}/4\\
+(\bar{b}\bar{\tau}(x)|\hat{\theta}|+2\hat{\theta}^{4}+2\hat{\theta}^{2}+2)\varpi{}^{2}.
\end{gather*}
 Let 
\begin{gather}
\alpha_{q}(x)=\rho^{4}(x)|x|^{6}+\bar{\pi}^{2}(x)\rho^{2}(x)|x|^{2} \nonumber\\
+\pi_{x}^{4}(x)+\bar{\pi}^{2}(x)m_{2}^{2}(x)x{}^{2}\nonumber\\
\alpha_{c}=\max\{\lambda^{2}/4,b^{2}/4,c^{2}/4,c^{4}/64\}\nonumber\\
\beta(z)=m_{1}^{4}(z)\|z\|^{4}+m_{1}^{2}(z)\|z\|^{2}\nonumber\\
\beta_{c}=\max\{1/4,\tau_{c}^{2}c^{2}/4\}\nonumber\\
\bar{\alpha}(x,\hat{\theta})=(\bar{b}\bar{\tau}(x)|\hat{\theta}|+2\hat{\theta}^{4}+2\hat{\theta}^{2}+2). \label{eq:bar_alpha}
\end{gather}
where $\alpha_c$ and $\beta_c$ are unknown.
As a result, one has
\[
\dot{U}(\varpi)\leq\alpha_{c}\alpha_{q}(|x|)+\beta_{c}\beta(\|z\|)+\bar{\alpha}(x,\hat{\theta})\varpi^{2}
\]
which further leads to 
\begin{equation}
\dot{U}(\varpi)\leq\bar{\alpha}_{\beta}\bar{\sigma}(\xi)\|\xi\|^{2}+\bar{\alpha}(x,\hat{\theta})\varpi{}^{2}\label{eq:Uw_dot_3}
\end{equation}
where $\bar{\alpha}_{\beta}=\max\{\alpha_{c},\beta_{c}\}$ is an unknown constant, and 
\begin{gather*}
\bar{\sigma}(\xi)=\rho^{4}(x)x{}^{4}+\bar{\pi}^{2}(x)\rho^{2}(x)+m_{3}^{4}(x)x^{2}\\
+\bar{\tau}^{2}(x)m_{2}^{2}(x)+m_{1}^{4}(z)\|z\|^{2}+m_{1}^{2}(z)
\end{gather*}

Let $\bar{V}(z,x,\tilde{\theta})=(1+1/b)V(z,x,\hat{\theta})$ where
$V(z,x,\hat{\theta})$ is given in (\ref{eq:Vzxt}) and its derivative
also satisfies $\dot{\bar{V}}(z,x,\tilde{\theta})\leq-\|\xi\|^{2}+\varpi^{2}/q.$\textcolor{red}{{}
}By Parameterized Changing Supply in Lemma \ref{lem:pcsf}, for any smooth
function $\Delta(\Xi)\geq\bar{\sigma}(\xi)+1$ and $\hat{k}=\max\{2\bar{\alpha}_{\beta},1\}$,
there exists a new supply function $\dot{V}_{\Xi}(\Xi)$ such that
\begin{gather}
\dot{V}_{\Xi}(z,x,\hat{\theta})\text{\ensuremath{\leq-\hat{k}\Delta(\Xi)\|\xi\|^{2}+\bar{p}'\varkappa(\tilde{\theta},\varpi)\varpi^{2}/q}}\nonumber\\
\text{\ensuremath{\leq-\hat{k}\Delta(\Xi)\|\xi\|^{2}+\bar{p}'\hat{\varkappa}(\theta)\bar{\varkappa}(\hat{\theta},\varpi)\varpi^{2}/q}} \label{eq:V_Xo}
\end{gather}
for some sufficiently smooth functions $\hat{\varkappa}$, $\bar{\varkappa}$
and unknown constant $\bar{p}'$. Note that $\hat{k}\Delta(\Xi)\geq2\bar{\alpha}_{\beta}\bar{\sigma}(\xi)+1$
and we specify $q$ to be $q\geq\bar{p}'\hat{\varkappa}(\theta)'$.
It leads to 
\begin{gather}
\dot{V}_{\Xi}(z,x,\hat{\theta})\text{\ensuremath{\leq-\hat{k}\Delta(\Xi)\|\xi\|^{2}+\bar{\varkappa}(\hat{\theta},\varpi)\varpi^{2}}}\label{eq:VXi}
\end{gather}

Now, it is ready to design the event-triggered law. Let 
\[
\bar{\gamma}(x,\hat{\theta},\varpi)=\bar{\varkappa}(\hat{\theta},\varpi)/2+\bar{\alpha}(x,\hat{\theta})+\frac{1}{2}
\]
where $\bar{\alpha}(x,\hat{\theta})$ is given in (\ref{eq:bar_alpha}) and $\bar{\varkappa}(\hat{\theta},\varpi)$ is given in (\ref{eq:V_Xo}).
One has 
\begin{gather*}
U(\varpi(t))=\int_{t_{k}^+}^{t}\dot{U}(\varpi(\tau))d\tau\leq\int_{t_{k}^+}^{t}\bar{\alpha}_{\beta}\bar{\sigma}(\xi(\tau))\|\xi(\tau)\|^{2}d\tau\\
+\int_{t_{k}^+}^{t}\bar{\alpha}(x(\tau),\hat{\theta}(\tau))\varpi{}^{2}(\tau)d\tau,\;t\in[t_{k},t_{k+1}].
\end{gather*}
which together with event-triggered law (\ref{eq:sd_con_ada-1}) (similar
to the proof of Theorem \ref{lem:uniform_lem}) implies 
\begin{gather*}
U(\varpi(t))\leq2\int_{t_{k}^+}^{t}\left(\bar{\alpha}_{\beta}\bar{\sigma}(\xi(\tau))\|\xi(\tau)\|^{2}d\tau\right.\\
\left.-\frac{1}{2}\bar{\varkappa}(\hat{\theta}(\tau),\varpi(\tau))\varpi^{2}(\tau)-\frac{1}{2}\varpi^{2}(\tau)\right)d\tau,\;t\in[t_{k},t_{k+1}].
\end{gather*}
Let $V(\Xi,\varpi)=V_{\Xi}(z,x,\hat{\theta})+U(\varpi)$ be the Lyapunov
function candidate for the closed-loop system (\ref{eq:gain_cl}).
As a result,
\begin{gather*}
V(\Xi(t),\varpi(t))-V(\Xi(t_k^+),\varpi(t_k^+))\leq\\
\leq-\int_{t_{k}^+}^{t}\|\xi(\tau)\|^{2}d\tau-\int_{t_{k}^+}^{t}\varpi^{2}(\tau)d\tau,\;\forall t\in[t_{k},t_{k+1}]
\end{gather*}
Therefore, $V(\Xi,\varpi)$ is monotonically decreasing and all signal
$z$, $x$, $\hat{\theta}$ and $\varpi$ are bounded. The convergence
and Zeno freeness can be proved similar to that of Theorem \ref{lem:uniform_lem}
and \ref{thm:robust} and are omitted here. \eproof

\section{Numerical Examples \label{sec:ne}}

\begin{figure}
\includegraphics[scale=0.45]{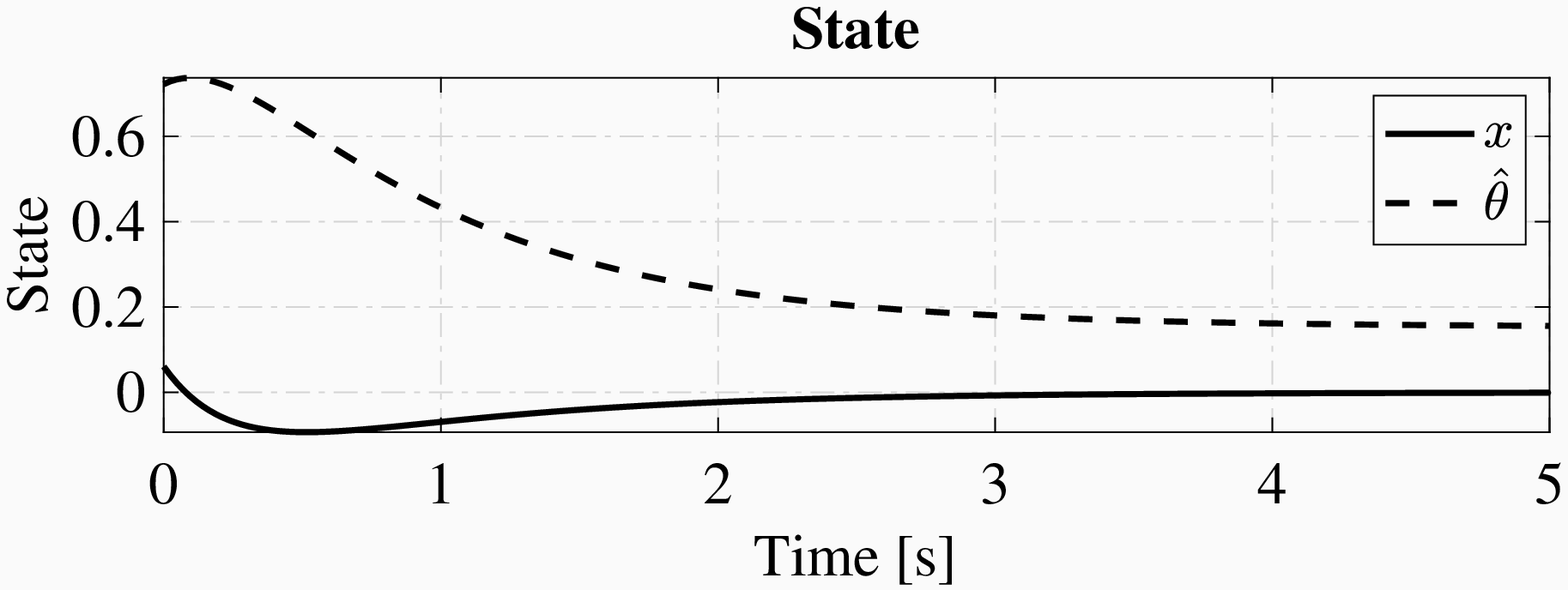}

\includegraphics[scale=0.45]{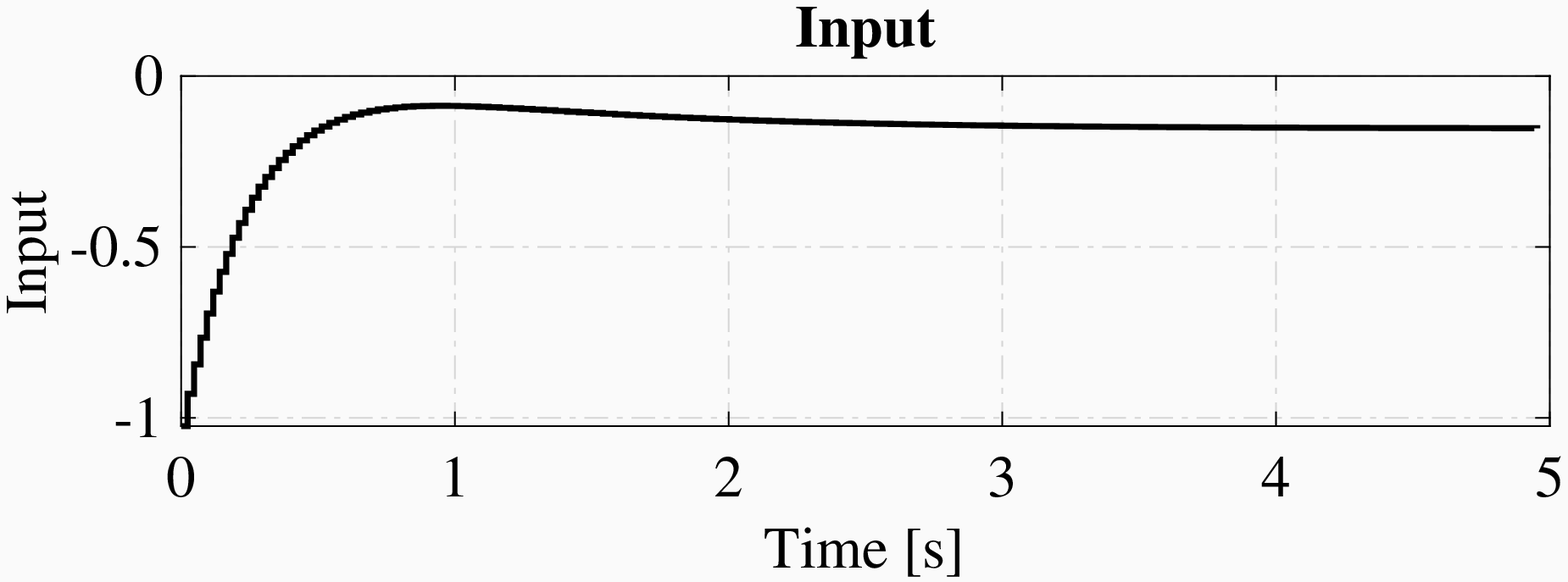}

\includegraphics[scale=0.45]{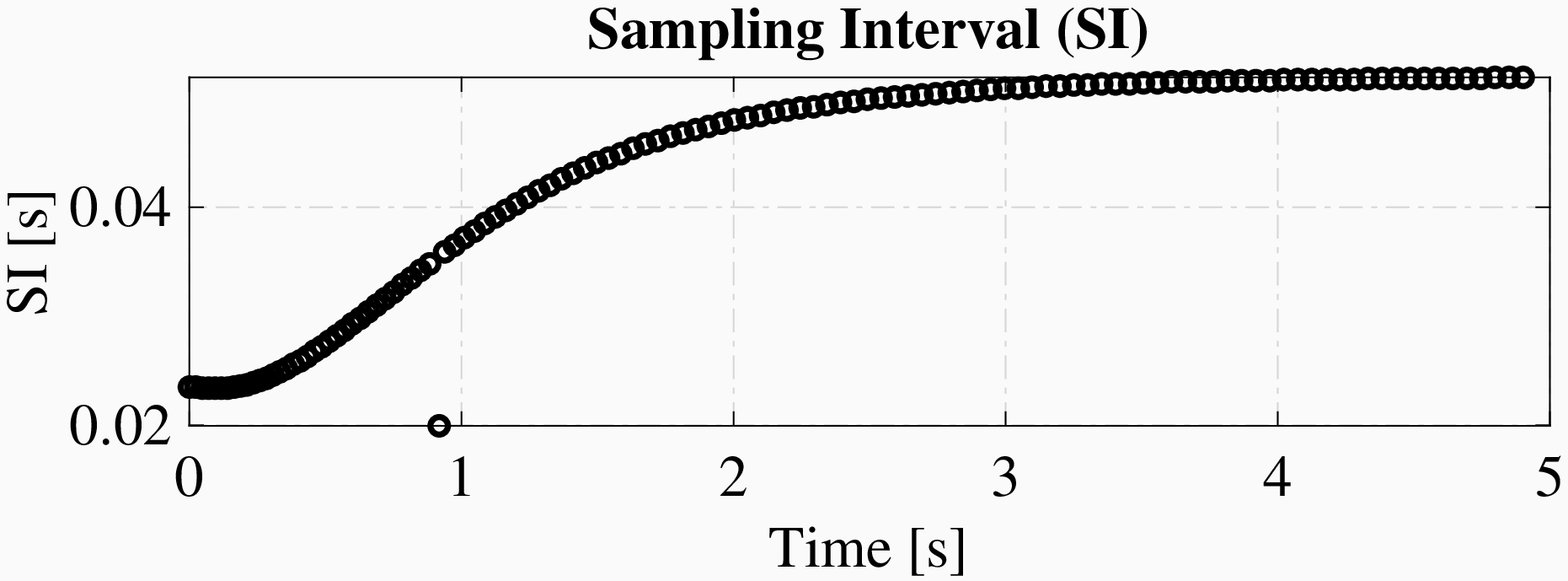}\caption{\label{fig:adaptive}The state trajectory, input profile and sampling interval of event-triggered
adaptive control in Example \ref{exp:1}. }

\end{figure}

\bexample\textcolor{red}{{} }\label{exp:1} (Event-triggered Adaptive
Stabilization with Uncertain Parameters) Consider the uncertain system
\[
\dot{x}=\theta\cos x+u
\]
where $\theta$ is an unknown parameter whose bounds is $|\theta|<2$.
We propose the sampling controller as follows,
\begin{eqnarray*}
u(t) & = & -\hat{\theta}(t_{k})\cos x(t_{k})-5x(t_{k})/4,\;t\in[t_{k},t_{k+1})\\
\dot{\hat{\theta}} & = & \lambda x\rho(x)+\lambda\varsigma(\hat{\theta},\varpi)
\end{eqnarray*}
with $\rho(x)$ and $\varsigma(\hat{\theta},\varpi)$ to be designed
and $\lambda>0$ to be specified. Let the sampling error $\varpi(t)$
be $\varpi(t)=-\cos(x(t_{k}))\hat{\theta}(t_{k})-5x(t_{k})/4+\cos(x(t))\hat{\theta}(t)+5x(t)/4,\;t\in[t_{k},t_{k+1})$.
Let $U(\varpi)=\frac{1}{2}\varpi^{2}$ and the inequality (\ref{eq:Uw_dot_ada})
is satisfied with $\dot{U}(\varpi)\leq\bar{\sigma}_{\varpi}(x)\|x\|^{2}+\bar{\alpha}_{\varpi}(x,\hat{\theta})\varpi{}^{2}-\lambda\varpi\varsigma(\hat{\theta},\varpi)\cos x+\frac{5}{4}\varpi\tilde{\theta}\cos x$
where $\bar{\sigma}_{\varpi}(x)=1.251+\frac{1}{4}\lambda^{2}$ and
$\bar{\alpha}_{\varpi}(x,\hat{\theta})=2\hat{\theta}^{4}+9\hat{\theta}{}^{2}+|\hat{\theta}|+\rho^{2}(x)-\frac{1}{4}$.
Let $V_{x}=\frac{1}{2}x^{2}$ and $\dot{V}_{x}\leq-x^{2}+\varpi^{2}-x\tilde{\theta}\cos x$.
Let $\Delta(x)=\lambda_{c}=3.501+\frac{1}{4}\lambda^{2}$ and thus
$\rho_{q}$ in (\ref{eq:rho_q}) should be $\rho_{q}\geq2\lambda_{c}$
and $\sigma(\varpi)=\rho_{q}$ Then, $\rho(x)$ and $\varsigma(\hat{\theta},\varpi)$
in (\ref{eq:ada_con}) is 
  $\rho(x)=\rho_{q}\cos x$ and 
$\varsigma(\hat{\theta},\varpi)=-\frac{5}{4}\varpi$.
As a result, $\tilde{\alpha}_{\varpi}(x,\hat{\theta})\geq\bar{\alpha}_{\varpi}(x,\hat{\theta})+\frac{5}{4}\lambda$
is selected to make (\ref{eq:alpha_w}) satisfied. Then, $\bar{\gamma}(x,\hat{\theta})$
can be selected as $\bar{\gamma}(x,\hat{\theta})=\tilde{\alpha}_{\varpi}(x,\hat{\theta})+\sigma(\varpi)/2+\frac{1}{2}=2\hat{\theta}^{4}+9\hat{\theta}{}^{2}+|\hat{\theta}|+\rho^{2}(x)+\frac{1}{4}+\frac{5}{4}\lambda+\frac{\rho_{q}}{2}$
with $\rho_{q}=2\lambda_{c}$. Finally, the event-triggered law is
designed according to Theorem \ref{thm:adap}, as follows,
\begin{gather*}
t_{k+1}=\inf_{t\geq t_{k}}\{(t-t_{k})\max_{\tau\in[t_{k},t]}\{\bar{\gamma}(x(\tau),\hat{\theta}(\tau))\varpi^{2}(\tau)\}\\
\geq \frac{1}{4}\|\varpi_{[t_{k},t]}\|^{2},\;\text{and }\|\varpi_{[t_{k},t]}\|\neq0\}.
\end{gather*}
The simulation result is illustrated in Fig \ref{fig:adaptive}, which
includes figures for state trajectories, input signal and sampling
intervals. Note that the estimated value $\hat{\theta}$ converges
to the real value of $\theta=0.1534$ and the trajectories of system
converge to the origin.

\begin{figure}
\includegraphics[scale=0.45]{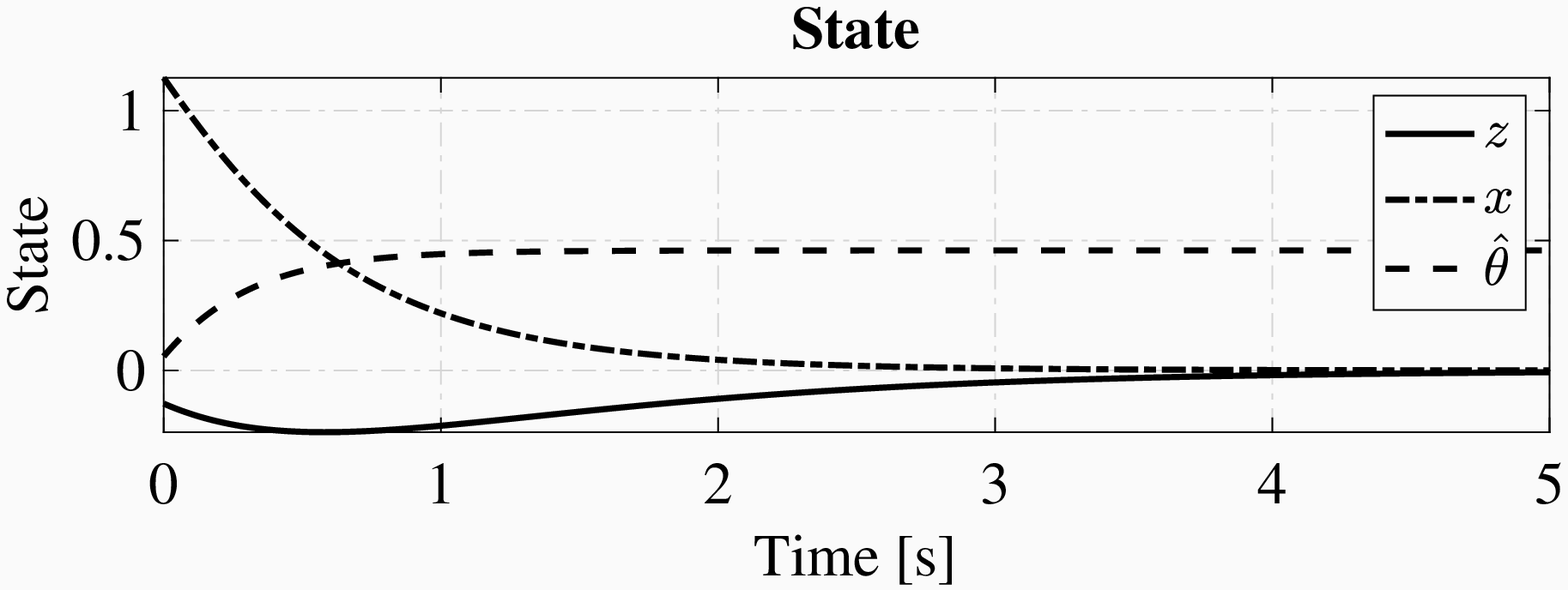}

\includegraphics[scale=0.45]{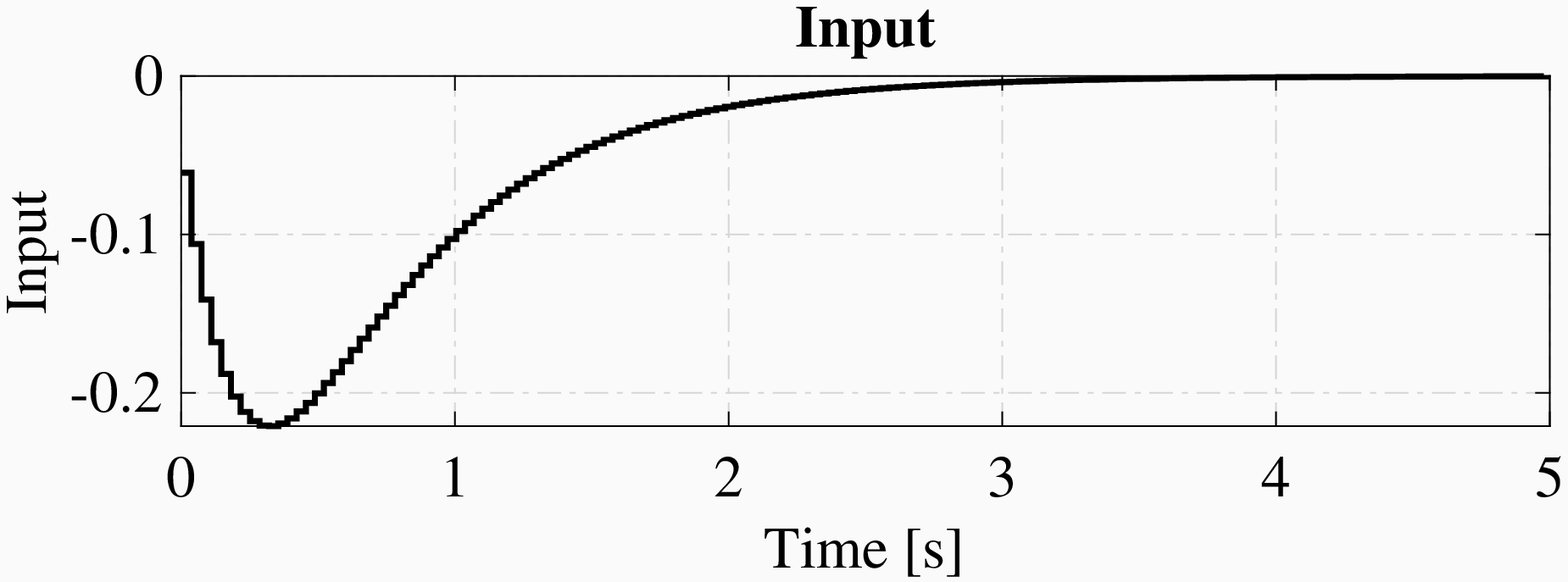}

\includegraphics[scale=0.45]{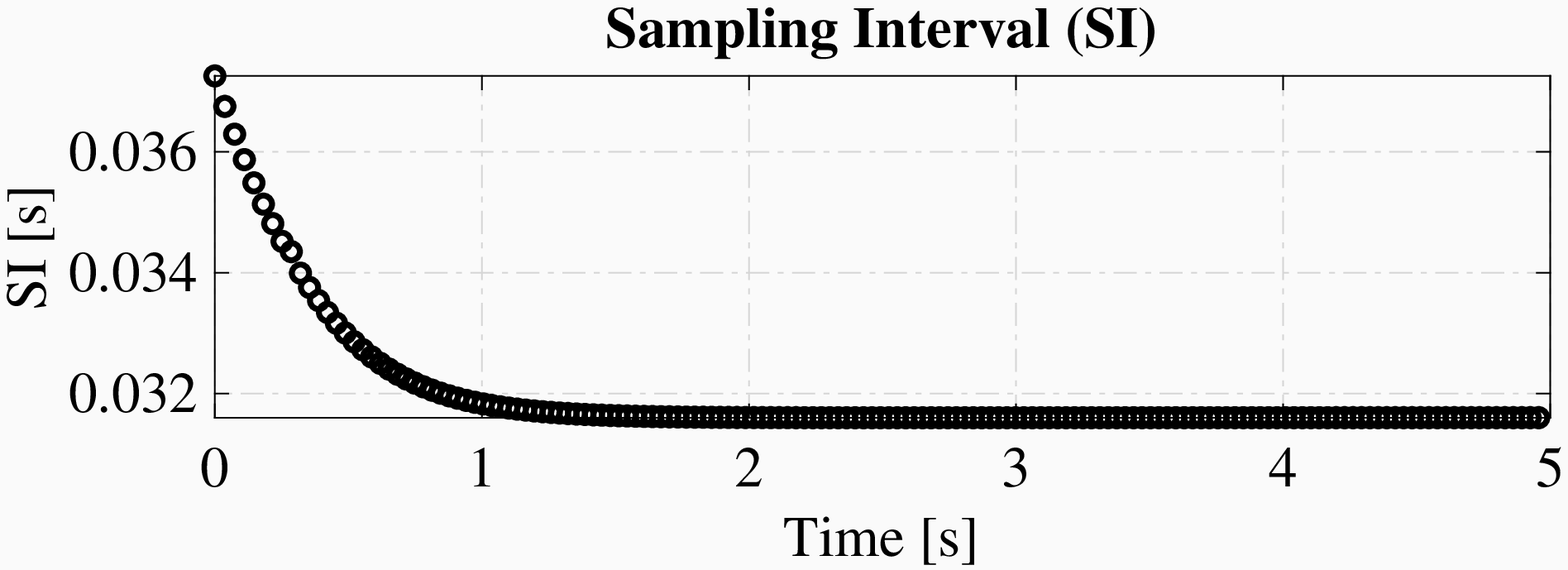}\caption{\label{fig:dg}The state trajectory, input profile and sampling interval of event-triggered
adaptive control in Example \ref{exp:2}.}
\end{figure}

\eexample \bexample \label{exp:2} (Event-triggered Stabilization
with Dynamic Gain) Consider the following nonlinear system 
\begin{eqnarray*}
\dot{z} & = & -z+w_{3}x\\
\dot{x} & = & w_{1}z\sin x+w_{2}x+bu
\end{eqnarray*}
where $w_{1}$, $w_{2}$ and $w_{3}$ are unknown parameters and $0<b<1$.
Note that when $V(z)=z^{2}$, one has $\dot{V}(z)\leq-z^{2}+px^{2}$
for any unknown constant $p\geq w_{3}^{2}$ and thus Assumption \ref{ass:za}
is satisfied. And $|w_{1}z\sin x+w_{2}x|\leq c(|z|+|x|)$ for some
$c\geq\max\{|w_{1}|,|w_{2}|\}$, that is, inequality (\ref{eq:f_bound_2})
is verified for $m_{1}=1$ and $m_{2}=1$. Then, we propose the sampling
controller as in (\ref{eq:dg_con}) and the closed-loop system is
written as 
\begin{eqnarray*}
\dot{z} & = & -z+w_{3}x\\
\dot{x} & = & w_{1}z\sin x+w_{2}x-b\hat{\theta}\rho(x)x+b\varpi,\\
\dot{\hat{\theta}} & = & \lambda\rho(x)x^{2},\;\lambda>0
\end{eqnarray*}
where $\varpi(t)=-\hat{\theta}(t_{k})\rho(x(t_{k}))x(t_{k})+\hat{\theta}(t)\rho(x(t))x(t)$
is the sampling error. Let $V(\Xi)=\frac{1}{2}x^{2}+\frac{1}{2}b\lambda^{-1}\tilde{\theta}^{2}+2z^{2}$
and $\dot{V}(z,x,\tilde{\theta})\leq-|z|^{2}+\left(c^{2}/4+c+b^{2}q/4+2p-b\theta\rho(x)\right)x^{2}+\varpi^{2}/q.$
Select $\rho(x)=1$ and $\theta\geq\left(c^{2}/4+c+2p\right)/b+bq/4.$
Then, $\dot{V}(z,x,\tilde{\theta})\leq-\|\xi\|^{2}+\varpi^{2}/q$.
Note that $\dot{\varpi}=-\dot{\hat{\theta}}x-2\hat{\theta}\dot{x}$
and let $U(\varpi)=\frac{1}{2}\varpi^{2}$. The inequality (\ref{eq:Uw_dot_3})
is satisfied, i.e., $\dot{U}(\varpi)\leq\bar{\alpha}_{\beta}\bar{\sigma}(\xi)\|\xi\|^{2}+\alpha(\hat{\theta})\varpi^{2}$
where
\begin{gather*}
\alpha(\hat{\theta})=1+\hat{\theta}^{2}+\hat{\theta}^{4}+2\bar{b}|\hat{\theta}|\\
\bar{\alpha}_{\beta}=\max\{\lambda^{2}/4,c^{2}+b^{2}\}\\
\bar{\sigma}(\xi)=1+\|\xi\|^{4}.
\end{gather*}
Note that $\frac{1}{2}x^{2}+\frac{1}{2}\lambda^{-1}\tilde{\theta}^{2}+2z^{2}\leq\bar{V}(z,x,\tilde{\theta})=(1+1/b)V(z,x,\tilde{\theta})\leq s\left[\frac{1}{2}x^{2}+\frac{1}{2}\lambda^{-1}\tilde{\theta}^{2}+2z^{2}\right]$
where $s\geq\max\{1+1/b,1+b\}$. Let $\hat{k}=\max\{2\bar{\alpha}_{\beta},1\}$
and $\Delta(\xi)=\bar{\sigma}(\xi)+1$. Let $\rho(\chi)=4+2\chi^{2}/\underline{\iota}^{2}$
and $V_{\Xi}(\Xi)=\int_{_{0}}^{\bar{V}(z,x,\tilde{\theta})}\hat{k}\rho(\tau)d\tau$.
Let $\bar{\iota}=\max\{\frac{1}{2}\lambda^{-1},2\}$ and $\underline{\iota}=\frac{1}{2}\min\{\lambda^{-1},1\}$.
The derivative of $V_{\Xi}(\Xi)$ becomes
\begin{gather*}
\dot{V}_{\Xi}(\Xi)\leq-\hat{k}\Delta(\xi)\|\xi\|^{2}+\frac{64\hat{k}s^{2}\bar{\iota}^{2}}{\underline{\iota}^{2}q}\left(|\theta|^{4}+1\right)\\
\times\left(\left(|\hat{\theta}|+|\varpi|\sqrt{2/q}\right)^{4}+1\right)
\end{gather*}
Let 
\[
q\geq\bar{c}\max\left\{ 2,\frac{64\hat{k}s^{2}\bar{\iota}^{2}}{\underline{\iota}^{2}}\left(|\theta|^{4}+1\right)\right\} .
\]
Then, inequality (\ref{eq:VXi}) is satisfied with $\bar{\varkappa}(\hat{\theta},\varpi)=\frac{\left(|\hat{\theta}|+|\varpi|\right)^{4}+1}{\bar{c}}$.
Then, $\bar{\gamma}(\hat{\theta},\varpi)=\bar{\varkappa}(\hat{\theta},\varpi)/2+\bar{\alpha}(\hat{\theta})+\frac{1}{2}$
where $\alpha(\hat{\theta})=1+\hat{\theta}^{2}+\hat{\theta}^{4}+2\bar{b}|\hat{\theta}|$.
Finally, the event-triggered law is designed according to Theorem
\ref{thm:dg}, as follows,
\begin{gather*}
t_{k+1}=\inf_{t\geq t_{k}}\{(t-t_{k})\max_{\tau\in[t_{k},t]}\{\bar{\gamma}(\hat{\theta}(\tau),\varpi(\tau))\varpi^{2}(\tau)\}\\
\geq\frac{1}{4}\|\varpi_{[t_{k},t]}\|^{2},\;\text{and }\|\varpi_{[t_{k},t]}\|\neq 0\}.
\end{gather*}
The simulation result is illustrated in Fig \ref{fig:dg}, which includes
figures of state trajectories, input signal and sampling intervals
showing that the system trajectories converge to the origin without Zeno behavior.

\eexample

\section{Conclusion \label{sec:con}}

In this paper, we propose a novel sampling control framework based
on the emulation technique where the sampling error is regarded as
the auxiliary input to the emulated system. The design of periodic
sampling and event-triggered control law utilizes the supremum norm
of sampling error and renders the error dynamics bounded-input-bounded-state
(BIBS), when coupled with system dynamics, achieves global or semi-global
stabilization. The proposed framework is then extended to tackle the
event-triggered and periodic sampling stabilization for systems with
dynamic uncertainties. The proposed framework is further extended
to solve two classes of event-triggered universal adaptive control
problems. It would be very interesting to consider the case where
the adaptation dynamics is also sampled and periodic event-triggered
control along this research line.

\appendix
\blem \label{lem:pcsf} (Parameterized Changing Supply Functions,
Lemma 6.1 in \cite{Chen2015Book}) Consider the system $\dot{\xi}=f(\xi,\varpi)$
with $\xi=\text{col}(z,x)\in\mathbb{R}^{n}$ and $\varpi\in\mathbb{R}^{m}$.
Suppose there exists a supplying function $V(\xi)$ satisfying 
\begin{gather*}
\underline{\alpha}(\|\xi\|)\leq V(\xi)\leq\bar{\alpha}(\|\xi\|,s)\\
\dot{V}(x)\leq-\alpha(\|x\|)+p\sigma(\|\varpi\|)
\end{gather*}
for $\mathcal{K}_{\infty}$ functions $\alpha$, $\underline{\alpha}$,
$\sigma$, a parameterized $\mathcal{K}_{\infty}$ function $\bar{\alpha}$,
and positive number $p$ and $s$. Then, for any smooth function $\Delta:\mathbb{R}^{n}\rightarrow[0,\infty)$
and positive number $k$, there exists a continuously differentiable
function $V'(\xi)$ satisfying 
\begin{gather*}
\underline{\alpha}'(\|\xi\|)\leq V'(\xi)\leq\bar{\alpha}'(\|\xi\|,s')\\
\dot{V}'(x)\leq-k\Delta(\xi)\alpha(\|x\|)+p'\varkappa(z,\varpi)\sigma(\|\varpi\|)
\end{gather*}
for a $\mathcal{K}_{\infty}$ functions $\underline{\alpha}'$, a
parameterized $\mathcal{K}_{\infty}$ function $\bar{\alpha}'$, a
$\mathcal{SN}$ function $\varkappa$ and positive numbers $p'$ and
$s'$. Moreover, if the functions $\alpha$, $\underline{\alpha},$
$\bar{\alpha}$, $\sigma$ and $\Delta$ are known, so are the functions
$\alpha'$, $\bar{\alpha}'$ and $\varkappa$. The positive numbers
$s$, $p$, $k$, $s'$ and $p'$ are not necessarily known. \elem

\prooflater{Theorem \ref{thm:psc}} It suffices to show condition
(\ref{eq:PSC1}) imply (\ref{eq:psc_con}) 
and meanwhile it is guaranteed that the boundedness of signals $\varpi$, i.e.,
$\|\varpi(t)\|\leq R_{0}$, $\forall t\geq t_{0}$. 
From (\ref{eq:PSC1}),
one has $\omega\underline{\alpha}_{\varpi}(s)-2\gamma(s)>0,\;\forall 0 <s <R_0$.
Since $\underline{\alpha}_{\varpi}$ and $\gamma$ are continuous
functions, there exists a sufficiently small $\tilde{R}$ such that
\[
\omega\underline{\alpha}_{\varpi}(R_{0}+s)-2\gamma(R_{0}+s)>0,\;\forall0\leq s\leq\tilde{R},
\]
It, together with (\ref{eq:PSC1}), further implies 
\begin{equation}
2\gamma(s)<\omega\underline{\alpha}_{\varpi}(s),\;\forall0<s\leq R_{0}+\tilde{R}.\label{eq:rq_ine}
\end{equation}

We will prove that $\|\varpi(t)\|\leq R_{0},\;\forall t\geq t_{0}$.
If this is not true, there exists a finite time $t_{a}>t_{0}$ such
that $\|\varpi(t)\|\leq R_{0}+\tilde{R},\;\forall t\in[t_{0},t_{a}]$
but $\|\varpi(t_{a})\|>R_{0}$. Due to (\ref{eq:rq_ine}) and the periodic sampling law in
(\ref{eq:period}), 
one  have 
\begin{gather}
2(t-t_{i})\gamma(\|\varpi\|_{[t_{i},t]})\leq 2T\gamma(\|\varpi\|_{[t_{i},t]}) <\underline{\alpha}_{\varpi}(\|\varpi\|_{[t_{i},t]})\nonumber \\ \leq\max_{\tau\in[t_{i},t]}\{U(\varpi(\tau))\},\;
\forall t\in[t_{i},t_{i}+T),\nonumber \label{eq:et_con}\\
\forall t_i=\{t_0,t_0+T,\cdots,t_0+(\bar{k}+1)T\},
\end{gather}
where 
 $\bar{k}=\max_{k}\left\{ k\mid t_{0}+(k+1)T\leq t_{a}\right\} $. Therefore,   (\ref{eq:psc_con})  
 is valid for $t< t_{a}$ and
 $T \leq \bar{t}_{k+1}-t_{k}$, $\forall k\in\mathbb{N}$
where $\bar{t}_{k+1}$ is given in (\ref{eq:t_k_bar}). By Remark \ref{rem:Tk},
results of Theorem \ref{lem:uniform_lem} hold, and one has  $V_{q}(x(t))\leq V_{q}(x(t_{0}))$. Due to $x(t_0)\in\mathcal{X}$, one has $V_{q}(x(t))\leq V_{q}(x(t_{0}))\leq \underline \alpha_q \circ \sigma_\varpi^{-1} \circ \alpha_\varpi(R_0)$, that is $\mathcal{S}\subseteq \mathbb{R}^n$ is a positively invariant set.  
By Remark \ref{rem:bound}, one has $\|\varpi(t)\|\leq R_{0},\forall t_{0}\leq t\leq t_{a}$
which causes a contradiction to $\|\varpi(t_{a})\|>R_{0}$. Therefore,
$\|\varpi(t)\|\leq R_{0},\;\forall t\geq t_{0}$. As a result, conducting  similar analysis  as above for $t\geq t_0$ shows that (\ref{eq:psc_con}) 
 is valid for $t\geq t_{0}$
and using the discussion in Remark \ref{rem:Tk}  
completes the proof. \eproof

\prooflater{Corollary \ref{cor:robust}} The fact that $\sigma_{\varpi}(s)=\mathcal{O}(\alpha(s))$
as $s\rightarrow0^{+}$ implies it is always possible to find a new
supply function $V_{q}(x)\sim\{\underline{\alpha}_{q},\bar{\alpha}_{q},\alpha_{q},(\sigma_{q},\varsigma_{q})\mid\dot{x}=f_{c}(x,\varpi,d)\}$
for some $\mathcal{K}_{\infty}$ functions $\underline{\alpha}_{q}$,
$\bar{\alpha}_{q}$, $\sigma_{q}$ and $\varsigma_{q}$, that can
be calculated accordingly. Denote $\upsilon=\text{col}(x,w)$. Following
the first proof step of Theorem \ref{lem:uniform_lem}, one has
\begin{gather*}
U(\varpi(t))\leq-2\int_{t_{k}^+}^{t}\hat{\alpha}_{\varpi}(\|\varpi(\tau)\|)d\tau+2\int_{t_{k}^+}^{t}\sigma_{\varpi}(\|x(\tau)\|)d\tau\\
+2\int_{t_{k}^+}^{t}\varsigma_{\varpi}(\|d(\tau)\|)d\tau,\;\forall t\in[t_{k},t_{k+1}].
\end{gather*}
Let $\bar{V}(x,\varpi)=V_{q}(x)+U(\varpi)$ be the Lyapunov function
candidate. Similarly, one has
\begin{gather*}
\bar{V}(x(t),\varpi(t))-\bar{V}(x(t_{k}^+),\varpi(t_{k}^+))\leq-\int_{t_{k}^+}^{t}\alpha_{\gamma}(\|x(\tau)\|)d\tau\\
-\int_{t_{k}^+}^{t}\gamma_{\alpha}(\|\varpi(\tau)\|)d\tau+\int_{t_{k}^+}^{t}\varsigma_{\upsilon}(\|d(\tau)\|)d\tau\\
\leq-\int_{t_{k}^+}^{t}\alpha_{\upsilon}(\|\upsilon(\tau)\|)d\tau+\int_{t_{k}^+}^{t}\varsigma_{\upsilon}(\|d(\tau)\|)d\tau,\;\forall t\in[t_k,t_{k+1}],
\end{gather*}
where functions $\alpha_{\upsilon}$ is selected as $\alpha_{\upsilon}(\|\upsilon\|)\leq\alpha_{\gamma}(\|x\|)+\gamma_{\alpha}(\|\varpi\|)$
and $\varsigma_{\upsilon}(s)\geq2\varsigma_{\varpi}(s)+\varsigma_{q}(s)$.
Note that there exist functions $\underline{\alpha}_{\upsilon}$ and
$\bar{\alpha}_{\upsilon}$ satisfying
\begin{gather*}
\underline{\alpha}_{\upsilon}(\|\upsilon\|)\leq\underline{\alpha}_{q}(\|x\|)+\underline{\alpha}_{\varpi}(\|\varpi\|),\\
\bar{\alpha}_{\upsilon}(\|\upsilon\|)\leq\bar{\alpha}_{q}(\|x\|)+\bar{\alpha}_{\varpi}(\|\varpi\|),
\end{gather*}
such that $\underline{\alpha}_{\upsilon}(\|\upsilon\|)\leq\bar{V}(x,\varpi)\leq\bar{\alpha}_{\upsilon}(\|\upsilon\|)$.
As a result, $\bar{V}(t)-\bar{V}(t_{k})\leq-\int_{t_{k}^+}^{t}\alpha_{\upsilon}\circ\bar{\alpha}_{\upsilon}^{-1}(\bar{V}(\tau))d\tau+\int_{t_{k}^+}^{t}\varsigma_{\upsilon}(\|d(\tau)\|)d\tau$
which implies that 
  any trajectory starting outside
$B$ will eventually goes inside $B$ where $B:=\{\upsilon\mid\bar{V}\leq\bar{\alpha}_{\upsilon}\circ\alpha_{\upsilon}{}^{-1}\circ\varsigma_{\upsilon}(\bar{d})\}$.
and $\|x\|\leq\underline{\alpha}_{\upsilon}^{-1}\circ\bar{\alpha}_{\upsilon}\circ\alpha_{\upsilon}{}^{-1}\circ\varsigma_{\upsilon}(\bar{d})$
is satisfied. Similar to the proof of Theorem \ref{lem:uniform_lem},
Zeno-freeness depends  on  the condition $\gamma(s)=\mathcal{O}(\underline{\alpha}_{\varpi}(s))$ as $s\rightarrow0^{+}$ and can be proved in the same way. \eproof

\bibliographystyle{plain}
\bibliography{lit}
\end{document}